%% file: a-vrs_JRSS.tex
\newtheorem{lemma}{Lemma}
\newtheorem{theorem}{Theorem}
\newtheorem{definition}{Definition}

\documentclass{statsoc}

\usepackage{graphicx}
\usepackage{url}
\usepackage{epstopdf}
\usepackage[usenames]{color}
\usepackage{rotating,afterpage}
\usepackage{multicol} 
\usepackage{amssymb}
\usepackage{latexsym}
\usepackage{amsmath}
\usepackage{natbib}
\def\argmax{\mathop{\rm arg\,max}} 

\title{An Auto-validating Rejection Sampler}
\author{Raazesh Sainudiin$^{\dag}$ and Thomas L.~York$^*$}
 \coaddress{current affiliation: Raazesh Sainudiin, 
	Department of Statistics,
	1 South parks Road, University of Oxford, 
	Oxford OX1 3TG, U.K.
	}
	\email{sainudii@stats.ox.ac.uk}
\address{$^\dag$Department of Mathematics and $^*$Department of Biological Statistics and Computational Biology, 
Cornell University, Ithaca, U.S.A.} 

\date{}	
\begin{document}

\maketitle

\begin{abstract}
In Bayesian statistical inference and computationally intensive frequentist inference, one is 
interested in obtaining samples from a high dimensional, and possibly multi-modal target density. 
The challenge is to obtain samples from this target without any knowledge of the normalizing 
constant.  Several approaches to this problem rely on Monte Carlo methods.  One of the simplest
such methods is the rejection sampler due to von Neumann.  Here we introduce 
an auto-validating version of the rejection sampler via interval analysis.  We show that our 
rejection sampler does provide us with independent samples from a large class of target densities 
in a guaranteed manner.  We illustrate the efficiency of the sampler by theory and by examples 
in up to 10 dimensions.  Our sampler is immune to the `pathologies' of some infamous densities 
including the witch's hat and can rigorously draw samples from piece-wise Euclidean spaces of small 
phylogenetic trees.

\end{abstract}

\input{a-vrs_core.tex}

\newpage

\bibliographystyle{plain}
\bibliography{references} 

\end{document}

%% file: a-vrs_core.tex
\section{INTRODUCTION}

Obtaining samples from a density $p(\theta) \triangleq p^*(\theta)/N_p$, 
where $\theta \in \mathbf{\Theta}$ and $\mathbf{\Theta}$ is a compact Euclidean subset, {\it i.e.,} 
$\mathbf{\Theta} \subset \mathbb{R}^n$, without any knowledge of the normalizing 
constant $N_p \triangleq \int_{\mathbf{\Theta}}{p^*(\theta)} \, d \theta$, is a basic problem in Bayesian 
inference and multivariate simulation.  Several approaches to 
this problem rely on computationally-intensive Monte Carlo methods through conventional 
floating-point arithmetic.  We will concentrate on the rejection sampler due to 
von Neumann \cite{Neumann1963}.  After a brief introduction to the rejection sampler (RS) in 
Section \ref{S:RSIntro}, an interval version of this sampler is formalized in Section \ref{S:MRS}.  This 
sampler is referred to as the Moore rejection sampler (MRS) in honor of Ramon E.~Moore who 
was one of the influential founders of interval analysis \cite{Moore1967}.  A brief introduction to 
interval analysis, a prerequisite to understanding MRS, as well as the notational conventions and 
background assumed in the rest of the paper, are given in Section \ref{S:AppA} for 
readers who are new to interval methods.  In Section \ref{S:AppB}, Lemma \ref{lemma1} 
shows that MRS produces independent samples from the desired target density and 
Lemma \ref{lemma2} describes the asymptotics of the acceptance probability for a refining 
family of MRSs.  Examples demonstrating the robustness and efficiency of MRS  
to complexity and dimensionality of the target are discussed in Section \ref{S:EX}.  We conclude in 
Section \ref{S:C}.  Our sampler is an adaptive and auto-validating von Neumann rejection sampler 
that can draw independent samples from a large class of target densities, including non-log-concave, 
sharp-peaked, and multi-modal targets.  Unlike many conventional samplers, each sample produced 
by MRS is equivalent to a computer-assisted proof that it is drawn from the desired target.
An open source {\tt C++} class library for MRS is publicly available from {\tt www.stats.ox.ac.uk/\url{~sainudii}/codes}.

\section{Rejection Sampler (RS)}\label{S:RSIntro}
Rejection sampling  \cite{Neumann1963} is a Monte Carlo method to draw independent 
samples from a target probability distribution $p(\theta) \triangleq p^*(\theta)/N_p$, 
where $\theta \in \mathbf{\Theta} \subset \mathbb{R}^n$.  Typically the target $p$ is any density that is 
absolutely continuous with respect to the Lebesgue measure.  In most cases of interest we 
can compute the target shape $p^*(\theta)$ for any $\theta \in \mathbf{\Theta}$, 
but the normalizing constant $N_p$ is unknown.  Given a proposal density $q=q^*/N_q$ and an envelope 
function $f_q$ (Figure \ref{F:RSMH}) over $\mathbf{\Theta}$ that satisfy certain conditions, RS can produce samples 
from $p$ as follows: 
\begin{figure}
\centering   \makebox[10pt]{
\input{pqfRSMH.tex}
}
\caption{The characteristics of three samplers with target $p = p^*/N_p$: (1) Rejection sampler with 
proposal $q=q^*/N_q$ and the envelope function $f_q$, (2) an independent 
Metropolis-Hastings sampler (IMHS) driven by an independent base chain $I$ 
with proposal $q_I = q^*_I / N_{q^*_I}$ and (3) a local Metropolis-Hastings sampler (LMHS) 
driven by a local base chain $L$ with proposal $q_L = q^*_L / N_{q^*_L}$   
centered at the current state (open square at the bottom). \label{F:RSMH}}
\end{figure}
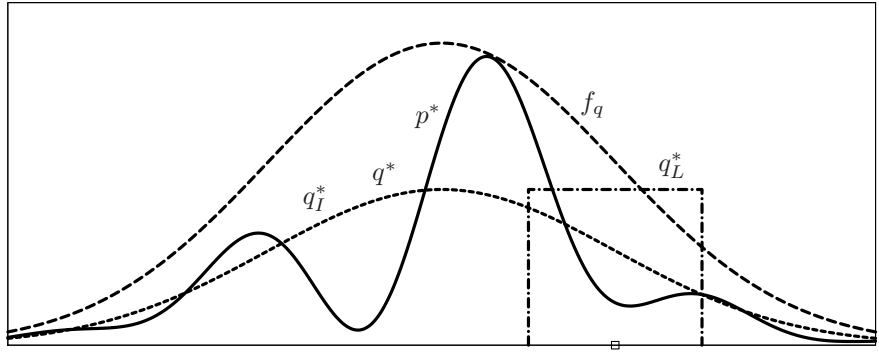

\subsection{Rejection Sampling Algorithm}\label{S:RS}
\begin{enumerate}
\item Choose a proposal density $q(\theta) = q^*(\theta)/N_q$ from which independent samples can be 
drawn, $N_q \triangleq \int_{\mathbf{\Theta}}{q^*(\theta)} \, d \theta$ is known, and $q^*(\theta)$ is computable for any $\theta \in \mathbf{\Theta}$.  

\item Find some $c$ for which the inequality
\begin{equation}
f_q(\theta) \triangleq c q^*(\theta) \geq p^*(\theta), \forall \, \theta \in \mathbf{\Theta} \label{E:EnvCond}
\end{equation}
is satisfied.  The smallest such value of $c$ is said to be optimal 
and denoted by $\hat{c}$, i.e.,
\[
\hat{c} \triangleq \inf \{c : c q^*(\theta) \geq p^*(\theta), \forall \, \theta \in \mathbf{\Theta}\}.
\]
\item {\label{Alg:RS}
Given (I)~a target shape $p^*(\theta)$, (II)~a proposal density $q(\theta)$, and (III)~an envelope 
function $f_q(\theta)$, $\forall \, \theta \in \mathbf{\Theta}$, 
that satisfy the above conditions, we can draw independent samples from the target $p(\theta)$ as follows:
\begin{enumerate}
\item GENERATE $T \sim q$. \label{RS1}
\item DRAW $H$ $\sim$ $Uniform[0,f_q(T)]$, where $f_q(T) \geq p^*(T)$. \label{RS2}
\item IF $H \leq p^*(T)$, THEN set $U = T$. 
\item RETURN to Step \ref{RS1}.
\end{enumerate}}
\end{enumerate}
It is not difficult to see that $U$ generated by the above algorithm is distributed 
according to $p$ \cite{Madras2002,Williams2001}.  
Observe that the probability $\mathbf{A}^p_{f_q}$ that a point proposed according to $q$ gets accepted 
as an independent sample from $p$ through the envelope function $f_q$ is the ratio of the integrals 
\[
\mathbf{A}^p_{f_q} = \frac{N_p}{N_{f_q}} \triangleq 
\frac{ \int_{\mathbf{\Theta}} {p^*(\theta) \, d \theta} }{\int_{\mathbf{\Theta}} {f_q(\theta) \, d \theta} },
\]
and the probability distribution over the number of samples from $q$ to obtain one sample from $p$ is 
geometrically distributed with mean $1/\mathbf{A}^p_{f_q}$ \cite{Madras2002,Williams2001}.  
Therefore, for a given $p$, we have to minimize $N_{f_q}$ over the allowed possibilities for $q$ and $f_q$ 
in order to obtain an efficient sampler with a high acceptance probability $\mathbf{A}^p_{f_q}$.

\section{Moore Rejection Sampler (MRS)}\label{S:MRS}

Moore rejection sampler (MRS) is an auto-validating rejection sampler (RS).  It can 
produce independent samples from any target shape $p^*$ that has a well-defined natural 
interval extension $P^*$ (Definition \ref{D:NIE}) over a compact domain $\mathbf{\Theta}$.  
MRS is said to be auto-validating because it automatically obtains a proposal $q$ that is easy 
to simulate from, and an envelope $f_q$ that is guaranteed to satisfy the envelope 
condition \eqref{E:EnvCond}.  MRS guarantees independent samples through auto-validating 
interval methods that also constitute the core of several recent computer-assisted proofs 
of challenging problems \cite{KulischBook2001}.  

\subsection{Theory}
In summary, the defining characteristics and notations of MRS are:
\[
\begin{array}{lcl}
\text{Compact domain} & \qquad & \mathbf{\Theta} = [\underline{\theta}, \overline{\theta}]\\
\text{Target shape} & \qquad & p^*(\theta) :\mathbf{\Theta} \rightarrow \mathbb{R}\\
\text{Target integral} & \qquad & N_p \triangleq \int_{\mathbf{\Theta}} {p^*(\theta) \, d \theta} \\ 
\text{Target density} & \qquad & p(\theta) \triangleq \frac{p^*(\theta)}{N_p} :\mathbf{\Theta} \rightarrow \mathbb{R}\\
\text{Interval extension of } p^* & \qquad & P^*(\Theta) : \mathbb{I}\mathbf{\Theta} \rightarrow \mathbb{IR}\\
\text{Proposal shape} & \qquad & q^*(\theta) :\mathbf{\Theta} \rightarrow \mathbb{R}\\
\text{Proposal integral} & \qquad & N_q \triangleq \int_{\mathbf{\Theta}} {q^*(\theta) \, d \theta} \\ 
\text{Proposal density} & \qquad & q(\theta) \triangleq \frac{q^*(\theta)}{N_q} :\mathbf{\Theta} \rightarrow \mathbb{R}\\
\text{Envelope function} & \qquad & f_q(\theta) = c q^*(\theta)\\
\text{Envelope integral} & \qquad & N_{f_q} \triangleq \int_{\mathbf{\Theta}} {f_q(\theta) \, d \theta}  = c N_q \\
\text{Acceptance probability} & \qquad & \mathbf{A}^p_{f_q} = \frac{N_p}{N_{f_q}}\\
\text{Partition of } \mathbf{\Theta} & \qquad & \mathfrak{T} := \{ \, \Theta^{(1)}, \Theta^{(2)}, ..., \Theta^{(|\mathfrak{T}|)}  \, \}.
\end{array}
\] 

If $p^* \in \mathfrak{E}$, the class of elementary functions (Definition \ref{D:ElemFunc}), 
and its natural interval extension $P^*$ is well-defined on $\mathbf{\Theta}$, then by Theorem \ref{3.1.11} 
\[
\begin{array}{lllllll}
Rng(p^* ; \mathbf{\Theta}) & \triangleq & p^*(\mathbf{\Theta}) & \subseteq & P^*(\mathbf{\Theta})  
& \triangleq &[\underline{P}^*(\mathbf{\Theta}), \overline{P}^*(\mathbf{\Theta})]
\end{array}
\]
which implies that
\begin{equation}
\underline{P}^*(\mathbf{\Theta}) \, \leq \, p^*(\theta) \, \leq \, \overline{P}^*(\mathbf{\Theta}), \, 
\forall \, \theta \in \mathbf{\Theta} \label{E:NaiveMRS}
\end{equation}
Although $[\underline{P}^*(\mathbf{\Theta}), \overline{P}^*(\mathbf{\Theta})]$ may 
over-estimate the range $p^*(\mathbf{\Theta})$, we can construct a naive 
MRS to draw samples from $p$ by using the following uniform proposal and constant envelope 
in Algorithm \ref{S:RS}.
\[
\begin{array}{lll}
q(\theta) & = & \frac{ \overline{P}^*(\mathbf{\Theta}) }{ d(\mathbf{\Theta}) \cdot \overline{P}^*(\mathbf{\Theta})} 
= \left( d(\mathbf{\Theta}) \right)^{-1}, \text{ and }\\
f_q(\theta) & = & \overline{P}^*(\mathbf{\Theta}),  
\end{array}
\]
where, $d(\mathbf{\Theta}) = d([\underline{\theta}, \overline{\theta}]) =  
\overline{\theta} - \underline{\theta}$ is the {\em diameter} of $\mathbf{\Theta}$. 
A lower bound for the acceptance probability of this naive MRS
is given by the range enclosure ratio: 
\[
\mathbf{A}^p_{\overline{P}^*(\mathbf{\Theta})}  =  \frac{N_p}{N_{\overline{P}^*(\mathbf{\Theta})}} = 
\frac{N_p}{d(\mathbf{\Theta}) \cdot \overline{P}^*(\mathbf{\Theta})} \geq 
\frac{d(\mathbf{\Theta}) \cdot \underline{P}^*(\mathbf{\Theta})}{d(\mathbf{\Theta}) \cdot \overline{P}^*(\mathbf{\Theta})} =
\frac{ \underline{P}^*(\mathbf{\Theta})}{\overline{P}^*(\mathbf{\Theta})}.
\]
Although this naive MRS can be extremely inefficient (i.e., can have a very low acceptance probability) for 
non-constant target shapes, one has the guarantee due to \eqref{E:NaiveMRS} that 
the necessary envelope condition \eqref{E:EnvCond} is satisfied.

A natural way to improve efficiency (i.e., increase the acceptance probability) is via partitions. 
Let $\mathfrak{T} \triangleq \{ \, \Theta^{(1)}, \Theta^{(2)}, ..., \Theta^{(|\mathfrak{T}|)}  \, \}$ be a finite 
partition of $\mathbf{\Theta}$.  Then by Theorem \ref{3.1.11} we can enclose $p^*(\Theta^{(i)})$, 
the range of $p^*$ over the $i$-th element of $\mathfrak{T}$, with the well-defined interval 
extension $P^*$ of $p^*$ over $\mathbf{\Theta}$
\begin{equation}\label{E:PartCont}
p^*(\Theta^{(i)}) \subseteq P^*(\Theta^{(i)}) \triangleq [\underline{P}^*(\Theta^{(i)}),\overline{P}^*(\Theta^{(i)})], \, \forall \, 
i \in \{1,2, ..., |\mathfrak{T}| \}. 
\end{equation}
For the given partition $\mathfrak{T}$ we can construct a partition-specific proposal $q^{\mathfrak{T}}(\theta)$ 
as a normalized simple function over $\mathbf{\Theta}$,
\begin{equation}\label{E:qRS}
q^{\mathfrak{T}}(\theta) = \left( N_{q^{\mathfrak{T}}} \right)^{-1} \, 
\sum_{i=1}^{|\mathfrak{T}|} { \overline{P}^*(\Theta^{(i)}) \, \mathbf{1}_{ \{ \theta \ \in \ \Theta^{(i)} \} } },
\end{equation}
where the normalizing constant is obtained from the sum 
\[
N_{q^{\mathfrak{T}}} \triangleq \sum_{i=1}^{|\mathfrak{T}|} \left( d(\Theta^{(i)}) \cdot \overline{P}^*(\Theta^{(i)}) \right).
\]  
The next ingredient $f_{q^{\mathfrak{T}}}(\theta)$ for our rejection sampler can simply be
\begin{equation}\label{E:fqRS}
f_{q^{\mathfrak{T}}}(\theta) = 
\sum_{i=1}^{|\mathfrak{T}|} { \overline{P}^*(\Theta^{(i)}) \, \mathbf{1}_{ \{ \theta \ \in \ \Theta^{(i)} \} } }
\end{equation}
The necessary envelope condition \eqref{E:EnvCond} is satisfied by $f_{q^{\mathfrak{T}}}(\theta)$
because of  \eqref{E:PartCont}.  Now, we have all the ingredients to perform a more efficient 
partition-specific Moore rejection sampling.  Lemma \ref{lemma1} shows that if the target 
shape $p^*$ has a well-defined natural interval extension $P^*$,
and if $U$ is generated according to the steps in part \ref{Alg:RS} of Algorithm \ref{S:RS}, 
and if the proposal density $q^{\mathfrak{T}}(\theta)$  and the envelope function 
$f_{q^{\mathfrak{T}}}(\theta)$ are given by  \eqref{E:qRS} and \eqref{E:fqRS}, respectively, 
then $U$ is distributed according to the target $p$.
Note that the above arguments as well as those in the proof of Lemma \ref{lemma1} naturally extend 
when $\mathbf{\Theta} \subset \mathbb{R}^n$ for $n>1$.  In the multivariate case, $\Theta^{(i)} \in \mathbb{IR}^n$ 
(Definition \ref{D:Boxes}) is a box.  Thus, we naturally replace the diameter of an interval by the {\em volume} of a box 
$v(\Theta^{(i)}) \triangleq  {\prod_{k=1}^{n}{d(\Theta_k^{i})} }$.  The envelopes and proposals are now simple
functions over a partition of the domain into boxes.  Analogous to the univariate case, the accepted samples 
are uniformly distributed in the region $S \subset \mathbb{R}^{n+1}$ `under' $p^*$ and `over' $\mathbf{\Theta}$.  
Hence their density is $p$ \cite{Williams2001}.

Next we bound the acceptance probability $\mathbf{A}^p_{f_{q^{\mathfrak{T}}}} \triangleq \mathbf{A}^p_{\mathfrak{T}}$ 
for this sampler.  Due to the linearity of the integral operator and \eqref{E:PartCont},
\[
\begin{array}{lcl}
N_p & \triangleq & \int_{\mathbf{\Theta}} {p^*(\theta) \, d \theta} \\
    &  = & \sum_{i=1}^{|\mathfrak{T}|} \int_{\Theta^{(i)}} {p^*(\theta) \, d \theta} \\
    &\in & \sum_{i=1}^{|\mathfrak{T}|} \left( d(\Theta^{(i)}) \cdot P^*(\Theta^{(i)}) \right) \\
    & =  & [ \ \, \sum_{i=1}^{|\mathfrak{T}|} \left( d(\Theta^{(i)}) \cdot \underline{P}^*(\Theta^{(i)}) \right) , \
    \sum_{i=1}^{|\mathfrak{T}|} \left( d(\Theta^{(i)}) \cdot \overline{P}^*(\Theta^{(i)}) \right) \, \ ].
\end{array}
\]
Therefore, 
\[
\mathbf{A}^p_{\mathfrak{T}} = \frac{N_p}{N_{f_{q^{\mathfrak{T}}}}} = 
\frac{N_p}{\sum_{i=1}^{|\mathfrak{T}|} \left( d(\Theta^{(i)}) \cdot \overline{P}^*(\Theta^{(i)}) \right)} \geq 
\frac{\sum_{i=1}^{|\mathfrak{T}|} \left( d(\Theta^{(i)}) \cdot \underline{P}^*(\Theta^{(i)}) \right)}
{\sum_{i=1}^{|\mathfrak{T}|} \left( d(\Theta^{(i)}) \cdot \overline{P}^*(\Theta^{(i)}) \right)}.
\]

We can say something more about the lower bound for $\mathbf{A}^p_{\mathfrak{T}}$ by 
limiting ourselves to target shapes within $\mathfrak{E_L}$, the Lipschitz class of elementary 
functions (Definition \ref{D:LipElemFunc}).  If $p^* \in \mathfrak{E_L}$ then we might expect the enclosure of $N_p$ to be proportional to the mesh $w$ of the partition $\mathfrak{T}$,
\[
w \triangleq \max_{i \in \{1,\dots,\mathfrak{T}\}}{d(\Theta^{(i)})}.
\]
Lemma \ref{lemma2} shows that if $p^* \in \mathfrak{E_L}$ and $\mathfrak{U}_W$ is a uniform partition 
of $\mathbf{\Theta}$ into $W$ intervals, then the acceptance probability 
$\mathbf{A}^p_{\mathfrak{U}_W} = 1 - \mathcal{O} (1/W)$.
More generally, any family of MRSs that construct their envelopes with  
$\overline{P}^*$ from the invoking family of refining partitions 
\[
\{ \ \mathfrak{T}_{\alpha} : \alpha \in \mathcal{A} \ \}
\]
can be thought of as a family of rejection samplers whose envelopes descend from 
above on $p^*$ in the form of simple functions.  The acceptance 
probability approaches $1$ at a rate that is no slower than linearly with the mesh.  
We can gain geometric insight into the sampler from an example.  The dashed lines of a given shade, 
depicting a simple function in Figure \ref{Fi:refine}, is a partition-specific envelope function \eqref{E:fqRS} 
for the target shape $s^{*}(x) = -\sum_{k=1}^5{ k \, x \, \sin{(\frac{k(x-3)}{3})}}$ over the domain $\mathbf{\Theta} = [-10,6]$ 
and its normalization gives the corresponding proposal function \eqref{E:qRS}.  As the refinement of $\mathbf{\Theta}$ 
proceeds through uniform bisections, the partition size increases as $2^i$, $i=1,2,3,4$.  Each of 
the corresponding envelope functions in increasing shades of gray can be used to draw auto-validated samples 
from the target $s(x)$ over $\mathbf{\Theta}$.  Note how the acceptance probability increases with refinement.

\subsection{Practice}
We theoretically studied the efficiency of uniform partitions for their tractability.  In practice, we may further 
increase the acceptance probability for a given partition size by adaptively 
partitioning $\mathbf{\Theta}$.  In our context, adaptive means the possible exploitation of any current 
information about the target.  We can refine the current partition $\mathfrak{T}_{\alpha}$ and 
obtain a finer partition $\mathfrak{T}_{\alpha^{\prime}}$ with an additional box
by bisecting a box $\Theta^{(*)} \in \mathfrak{T}_{\alpha}$ along the side with the maximal diameter.  
There are several ways to choose a $\Theta^{(*)} \in \mathfrak{T}_{\alpha}$ for bisection.  We 
explore three ways of choosing $\Theta^{(*)}$ from the current partition: (a) the box with the
largest volume, (b) the box with the largest diameter for its range enclosure and (c) the box
with the largest diameter for the product of its volume and its range enclosure.  When 
$\Theta^{(i)} \in \mathbb{IR}^n$ with volume 
$v(\Theta^{(i)})$, the three schemes can be 
formalized as follows:
\begin{alignat}{3} \label{E:PartScheme}
(a) &\text{ Volume-based} \qquad && \Theta^{(*)} = \argmax_{\Theta^{(i)} \in \mathfrak{T}_{\alpha} } 
{v(\Theta^{(i)})} \notag \\
(b) &\text{ Range-based} \qquad &&  \Theta^{(*)} =  \argmax_{\Theta^{(i)} \in \mathfrak{T}_{\alpha}} 
{ d(P^*(\Theta^{(i)}) ) }  \\
(c) &\text{ Integral-based}  \qquad &&  \Theta^{(*)} =  \argmax_{\Theta^{(i)} \in \mathfrak{T}_{\alpha}} 
{  \left( v(\Theta^{(i)}) \cdot d(P^*(\Theta^{(i)}) \right) } \notag
 \end{alignat}

Given a partitioning scheme, we employ a priority queue to conduct sequential refinements 
of $\mathbf{\Theta}$.  This approach avoids the exhaustive $\argmax$ computations to obtain
the $\Theta^{(*)}$ for bisection at each refinement step.  
A priority queue (PQ) is a container in which the elements may have different user-specified priorities.  
The priority is based on some sorting criterion that is applicable to the elements in the container.  
The PQ can be thought of as a collection in which the ``next'' element is always the
one with the highest priority, i.e., the largest with respect to the specified sorting criterion.  Since 
this container sorts using a {\em heap} which can be thought of as a binary tree, one can add 
or remove elements in logarithmic time.  This is a desirable feature of the PQ.  We implement the 
above three refinement schemes 
through PQs based on their respective sorting criterion.  The (a) the volume-based PQ 
manages the family of partitions $\mathfrak{U}_{W}$, (b) the range-based PQ manages 
the family $\mathfrak{R}_{\alpha}$ and (c) the integral-based PQ manages the family 
$\mathfrak{V}_{\alpha}$.

Once we have any partition $\mathfrak{T}$ of $\mathbf{\Theta}$, we can efficiently sample 
$\theta \sim q^{\mathfrak{T}}$ given by \eqref{E:qRS} in two steps.  First we sample a box 
$\Theta^{(i)} \in \mathfrak{T}$ according to the discrete distribution $t(\Theta^{(i)})$,
\begin{equation}
t(\Theta^{(i)}) =
 \frac{v(\Theta^{(i)}) \cdot \overline{P}^*(\Theta^{(i)})}
 {\sum_{i=1}^{|\mathfrak{T}|}{v(\Theta^{(i)}) \cdot \overline{P}^*(\Theta^{(i)})}},
  \   \Theta^{(i)} \in \mathfrak{T}, 
\end{equation}
and then we choose a $\theta \in \Theta^{(i)}$ uniformly at random.  Sampling from large discrete 
distributions (with million states or more) can be made faster by preprocessing the probabilities 
and saving the result in some convenient lookup table.  This basic idea \cite{Marsaglia1963} 
allows samples to be drawn rapidly.  We employ a more efficient preprocessing strategy \cite{Walker1977} 
that allows samples to be drawn in constant time even for very large discrete distributions as 
implemented in the GNU Scientific Library \cite{GSL}.  Thus, by means of priority queues and lookup
tables we can efficiently manage our adaptive partitioning of the domain for envelope construction, 
and rapidly draw samples from the proposal distribution.  We used the Mersenne Twister 
random number generator \cite{Matsumoto1998} in this paper.  Our sampler class builds 
on {\tt C-XSC 2.0}, a {\tt C++} class library for extended scientific computing using interval 
methods \cite{Hofschuster2004}.  All computations were done on a 2.8 GHz Pentium IV machine with
1GB RAM.  Having given theoretical and practical considerations to our 
Moore rejection sampler, we are ready to draw samples from various targets. 

\section{Discussion with Examples}\label{S:EX}

We empirically study sampler efficiency by sampling from qualitatively diverse targets since 
analytical results on efficiency are sharp only for relatively simple target parameterizations. 
In Section \ref{S:UGM} we first study the relative efficiencies of MRSs managed by the three 
PQs \eqref{E:PartScheme} by sampling from univariate Gaussian mixture targets.  Next, we study
the effects of target complexity (number of components, scales and domain size) 
on sampler efficiency.  In Section \ref{S:BL} we study the sampler behavior for a highly multi-modal 
two-dimensional target that is sensitive to a temperature parameter.  Using a trivariate mixture target 
in Section \ref{S:TNH}, we compare MRS to Monte Carlo Markov chain (MCMC) methods that rely 
on heuristic convergence diagnostics and exploit the connections between RS, importance sampler 
(IS) and independent Metropolis-Hastings sampler (IMHS) to simultaneously produce samples from 
all of them.  The effect of dimensionality on sampler efficiency is studied in Section \ref{S:MRWH}
where we draw samples from multivariate targets, including the multivariate witch's hat.  Section \ref{S:JCTQ} 
extends the sampler to piece-wise Euclidean domains of tree spaces.  Here we draw auto-validating 
samples of small trees (triplets and quartets) from the target likelihood function based on primate 
molecular sequence data.

\subsection{Univariate Gaussian Mixture}\label{S:UGM}
We apply MRS to targets whose shape $g_n$ is obtained from finite mixtures of
$n$ univariate Gaussian densities truncated over an interval $\mathbf{\Theta}$.  The means 
($\mu_i$'s), standard deviations ($\sigma_i$'s), weights ($w_i$'s), and domains ($\mathbf{\Theta}$'s) 
for each of the six targets studied are shown in Table \ref{T:3M-R.Exs}.  

\begin{table}
\caption{\label{T:3M-R.Exs} Moore rejection sampling from six different Gaussian mixture target 
shapes $g_n$ truncated over $\mathbf{\Theta}$, 
where $n$ is the number of mixture components.}
\centering
\begin{tabular}{|c c l|}
\hline
Target & $\mathbf{\Theta}$ & Parameters \\ \hline 
$g_1(x)$ & $[-10^2,10^2]$ & $\mu_1=-5,\sigma_1=1$, and $w_1=1.00$ \\ 
$g_2(x)$ & $ [-10^2,10^2]$ & $\mu_1=-5,\sigma_1=1$, $w_1=0.25$, $\mu_2=50$, \\
	&		& $\sigma_2=0.25$ \\ 
$g_5(x)$ & $[-10^2,10^2]$ & $\mu_1=-15, \mu_2=-5, \mu_3=3, \mu_4=6, \mu_5=50$, \\
	&		&	$\sigma_1=\sigma_2 = \sigma_4=1$, $\sigma_3=0.5$, $\sigma_5=0.1$,\\
	&		&	 $w_1=0.15, w_2=0.2, w_3=0.05, w_4=0.1$ \\ 
$g_5^{\prime}(x)$ & $[-10^2,10^2]$ & same as $g_5(x)$, except \\
		&		& $\sigma_1=\sigma_2 = \sigma_4=0.1$, $\sigma_3=0.05, \sigma_5=0.01$ \\
$g_5^{\prime \prime}(x)$ & $[-10^2,10^2]$ & same as $g_5(x)$, except  $\sigma_1=\sigma_2 = \sigma_4=0.01$,\\ 
& & $\sigma_3=0.005, \sigma_5=0.001$  \\
$\widehat{g}_5(x)$ & $[-10^{100},10^{100}]$ & same as $g_5(x)$ \\ \hline
\end{tabular}
\end{table}
 
\begin{figure}
\input{Gauss5_compare_PSzVar.tex}
\caption{Acceptance probability ($\mathbf{A}^{p}_{\mathfrak{T}_{\alpha}}$) versus partition 
size ($|\mathfrak{T}_{\alpha}|$) for six target shapes 
$p^* = g_1, g_2, g_5, g_5^{\prime}, g_5^{\prime \prime}, \widehat{g}_5$ (Table \ref{T:3M-R.Exs}) 
under different families of partitions: (1) volume-based $\mathfrak{U}_{W}$, (2) range-based 
$\mathfrak{R}_{\alpha}$ and (3) integral-based $\mathfrak{V}_{\alpha}$ 
(see text for description).\label{F:Gauss5CompM-R}}
\end{figure}
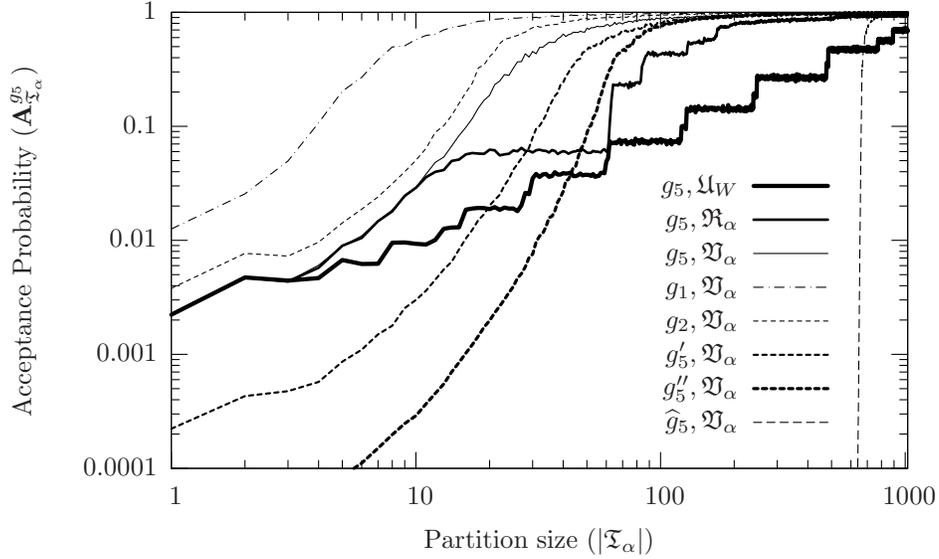

First, we study the efficiency of the three partitioning schemes \eqref{E:PartScheme} by
Moore rejection sampling from $g_5$.  Figure \ref{F:Gauss5CompM-R} shows the empirical 
acceptance probability of MRS, calculated from up to $10,000$ draws from a maximum of 
$100,000$ trials,  at each  partition size $|\mathfrak{T}_{\alpha}|$ for 
each of the three different families of partitions ($\mathfrak{U}_{W}$, $\mathfrak{R}_{\alpha}$ 
and $\mathfrak{V}_{\alpha}$).  Thus, for a given partition size $|\mathfrak{T}_{\alpha}|$, the domain interval 
$\mathbf{\Theta}$ gets adaptively partitioned through $|\mathfrak{T}_{\alpha}|-1$ 
bisections by the appropriate PQ.  The family of partitions $\mathfrak{V}_{\alpha}$ managed by 
the integral-based PQ is the most efficient as it can direct the next refining bisection towards 
the interval with the most uncertainty in its integral estimate.  The efficiency of the integral-based 
scheme is even more pronounced for multivariate exponential mixtures (results not shown).

Note that Lemmas \ref{lemma1} and \ref{lemma2} 
guarantee that MRS produces independent draws from any target in $\mathfrak{E_L}$.  
This includes Gaussian mixture targets with any finite number of components truncated over any 
compact interval.  Furthermore, the locations inside $\mathbf{\Theta}$ are arbitrary and the 
scales can be highly spiked (i.e., provided $\sigma_i >  0$ and can be enclosed by 
a machine interval via directed rounding \cite{Hammer1995,Kulisch2001}).  However, the efficiency
of the sampler can depend on (i) number of components, (ii) spikiness of peaks and (iii) domain size.  
We empirically study these effects by sampling from the six targets 
(Table \ref{T:3M-R.Exs}) using the family of MRSs induced by the most efficient 
partitions $\mathfrak{V}_{\alpha}$.  The acceptance probability plots (Figure \ref{F:Gauss5CompM-R}) 
for targets $g_1$, $g_2$, and $g_5$ illustrate the diminishing effect of the number of components on efficiency and
those for targets $g_5$, $g_5^{\prime}$ and $g_5^{\prime\prime}$, with progressively smaller variances, 
illustrate a similar effect of spikiness on efficiency at every partition size 
$|\mathfrak{V}_{\alpha}|$.  Note that in both cases 
sampler efficiency quickly recovers for larger partition sizes ($>100$).  Next we study the effect of domain size.  
In a computer, we cannot represent 
the real line and are forced to approximate it with the entire number screen, a compact interval.  Thus, the domain
of any target is necessarily truncated in a machine.  The acceptance probability plot for the target 
shape $\widehat{g}_5$, that is obtained by extending the domain of $g_5$ to a large interval of 
radius $10^{100}$ centered at $0$, shows the effect of domain size.  The first $700$ bisections or so 
are spent on zoning in on the intervals with relatively higher probability mass.  However, 
by $1000$ bisections our acceptance probability is almost $1$.

\subsection{Bivariate Levy}\label{S:BL}

\begin{figure}
\makebox{\centerline{\includegraphics[width=6.0in]{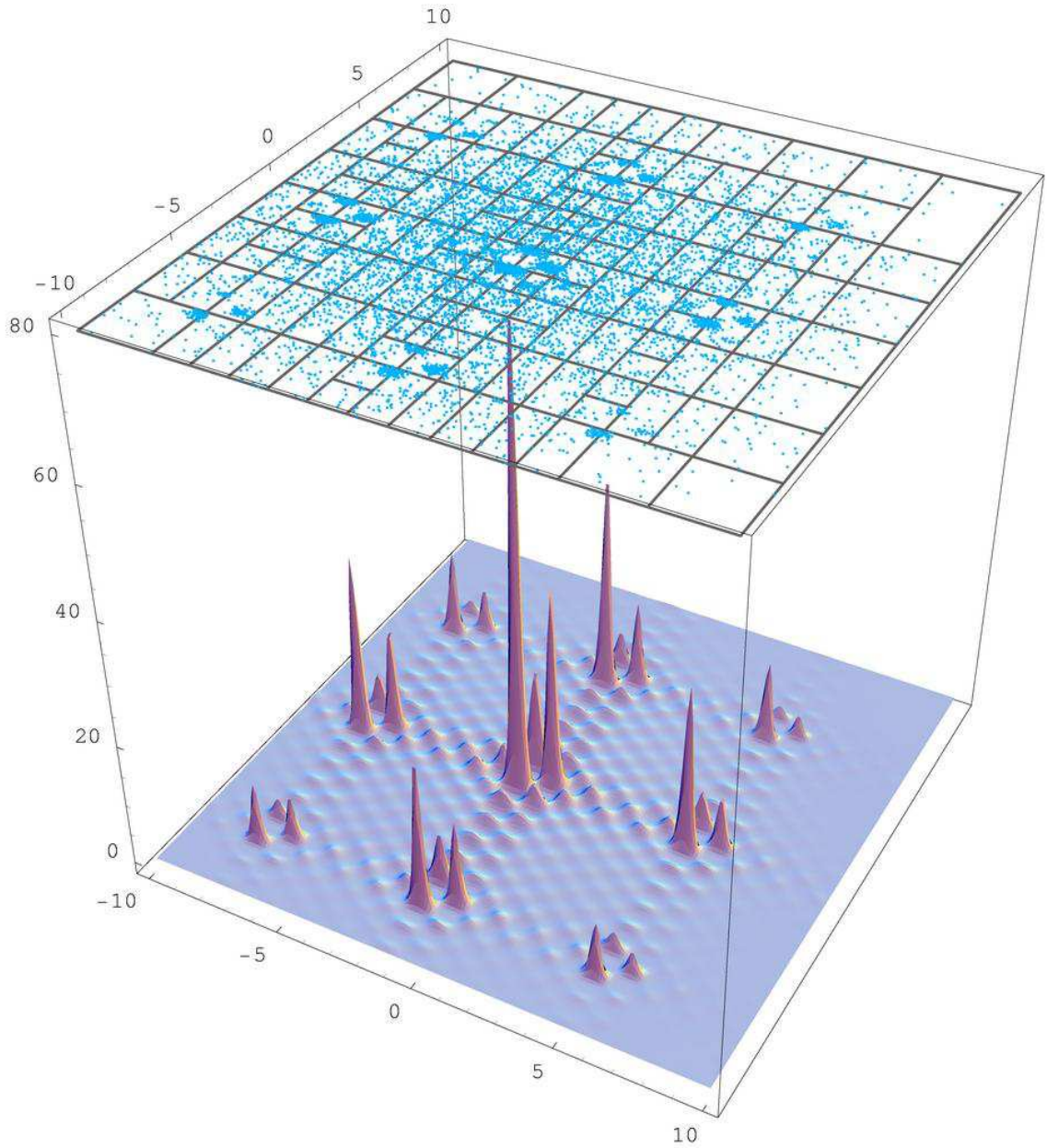}}} 
\caption{Shape of the Levy density $l^*_{40}$ with its 
$700$ modes \eqref{E:LevyShape}.  $10,000$ samples (points on top) from $l_{40}$ using the MRS 
induced by an adaptive partitioning of the domain into $150$ rectangles (with gray boundaries).
\label{F:LevyParts}}
\end{figure}

The bivariate Levy density $l_T (X_1,X_2)$ over 
$\mathbf{\Theta} \triangleq (\mathbf{\Theta}_1,\mathbf{\Theta}_2)^{\prime} = [-100,100]^{\otimes 2}$ 
\eqref{E:LevyShape} with temperature parameter $T$ and normalizing constant $N_{l_T}$ has 700 modes.
Figure \ref{F:LevyParts} shows $l^*_{40}$, i.e., the shape of the Levy density when $T=40$
and $10,000$ samples drawn from $l_{40}$ using the MRS induced by an integral-based adaptive 
partitioning of the domain into $150$ rectangles.   This MRS produced $10,000$ independent samples 
in less than $10$ CPU seconds at an acceptance probability of about $0.01$.  Mixtures of bivariate Gaussian 
shapes yielded comparable results.
\begin{alignat}{2} \label{E:LevyShape}
l_T (X_1,X_2) 	 =  &\frac{1}{N_{l_T}} l^*_T, \ \text{where, } \ l^*_T =  \exp \{ {-E(X_1,X_2)/T} \}, \\
E(X_1,X_2)     	 =  & \sum_{i=1}^5 {i \cos{((i-1)X_1+i)} } \sum_{j=1}^5 {j \cos{((j+1)X_2+j)} } \notag \\
			     & 
			     + (X_1 + 1.42513 )^2 + (X_2 + 0.80032)^2. \notag
\end{alignat}

\begin{figure}
\input{Levy_compare_Temp.tex}
\caption{Acceptance probability ($\mathbf{A}^{l_T}_{\mathfrak{V}_{\alpha}}$) and CPU seconds 
versus partition size ($|\mathfrak{V}_{\alpha}|$) for Levy targets $l_T$, where $T$ is 
the temperature \eqref{E:LevyShape}. \label{F:LevyComp}}
\end{figure}
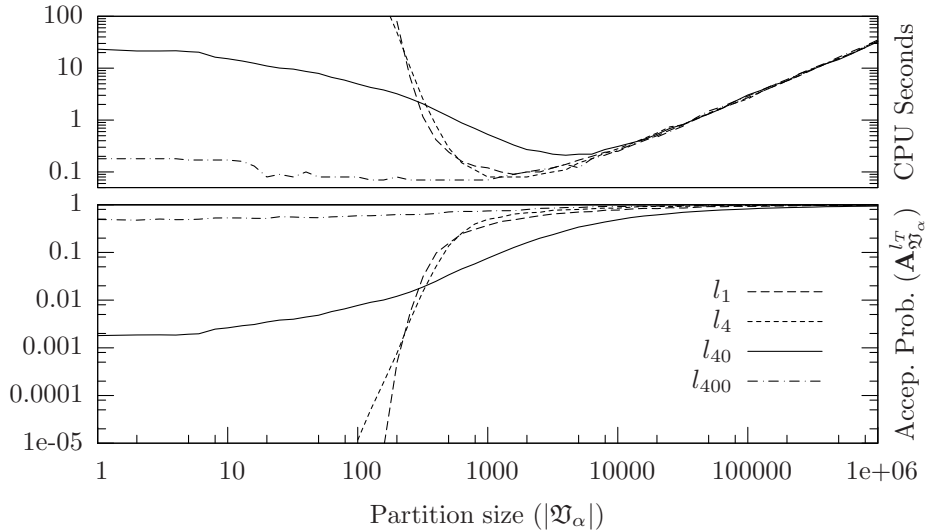

As the temperature parameter $T$ in $l_T$ increases, the density approaches a uniform
distribution on $\mathbf{\Theta}$.  The density is more peaked at low values of $T$.  Various
MCMC methods that use local proposals tend to mix well at higher temperatures and get
trapped at local peaks when $T$ is small.  To study the effect of temperature on our sampler's
efficiency, we plot the empirical acceptance probability as well as the CPU seconds taken 
to draw $10,000$ samples from each of four Levy targets at different temperatures 
($T=1,4,40,400$) as a function of the partition size $|\mathfrak{V}_{\alpha}|$ (Figure \ref{F:LevyComp}).  
The efficiency decreases as the temperature cools.  However, across the range of $T$ we explored, MRS 
can produce $10,000$ independent samples from $l_T$ in a guaranteed manner within $10$ CPU 
seconds with an acceptance probability greater than $1/100$.  Note that it is difficult to get a 
Monte Carlo Markov chain to mix properly and even more difficult to prove convergence
for such targets.

\subsection{Trivariate Needle in the Haystack}\label{S:TNH}

Using the target shape $h^*$  \eqref{E:Needle} over $\mathbf{\Theta} = [-10,10]^{\otimes 3}$, we compare
MRS to a popular MCMC sampler that relies on heuristics for convergence 
diagnosis and exploit the connection between three Monte Carlo methods.  We begin with an 
introduction to an MCMC sampler known as the 
Metropolis-Hastings sampler (MHS) \cite{Metropolis1953,Hastings1970} 
and a commonly used statistic for convergence diagnosis.
\begin{equation} \label{E:Needle}
h^*(x) = \frac{1}{\sigma_1^3} \exp \{ {-\frac{1}{2} ((x - \mu_1)/
\sigma_1)^2} \} + 
         \frac{1}{\sigma_2^3} \exp \{ {-\frac{1}{2} ((x - \mu_2)/
\sigma_2)^2} \}
\end{equation}

\subsubsection{ Metropolis-Hastings Sampler and Convergence Diagnostics}

Given $q_Y(\theta,\cdot)$, a possibly dependent proposal distribution for the base Markov chain $Y$
(Figure \ref{F:RSMH}), the following algorithm produces a Markov chain known as the 
Metropolis-Hastings (MH) chain on $\mathbf{\Theta}$ merely from the knowledge of ratios of
the form $p^*(\theta)/p^*(\theta^{\prime})$ for any 
$(\theta,\theta^{\prime}) \in \mathbf{\Theta} \times \mathbf{\Theta}$.  The stationary distribution of 
the MH chain is $p$.
\begin{itemize}
\item[1] Choose an arbitrary starting point $\theta_0$ and set $i=0$.
\item[2] Generate a candidate point $\theta^{\prime}$ $\sim$ $q_Y(\theta_i,\cdot)$ and $u$ $\sim$ $U(0,1)$.
\item[3] {Set:
\begin{equation}
\theta_{i+1} = 
\begin{cases}
{\theta^{\prime}} & \text{if } u \leq \, \frac{p^*(\theta^{\prime}) q_Y(\theta^{\prime},\theta_i)}{p^*(\theta_i) q_Y(\theta_i,\theta^{\prime})} \notag \\
\theta_i  & \text{otherwise} \notag
\end{cases}
\end{equation}}
\item[4] Set $i=i+1$ and GO TO 2
\end{itemize}
When the base chain has an independent proposal we refer to the MH chain it drives according to
the above algorithm as the independent Metropolis-Hastings sampler (IMHS) and when the base
chain has a local proposal we refer to the corresponding MH chain as the local Metropolis-Hastings 
sampler (LMHS).

Although the MH chain asymptotically approaches $p$, it is not trivial to know if it has converged even for 
relatively simple cases \cite{Diaconis1998}.  
One often resorts to some heuristic convergence diagnostics.  A fairly popular diagnostic statistic
\cite{Gelman1992,Gelman1996} runs multiple MH chains with randomly dispersed initial 
conditions and compares the within ($W$) and between ($B$) chain variation of the sampled draws.  When the 
ratio $B/W$ is small enough, one can be fairly certain that all the chains have converged to the same 
distribution.  Note that $B/W$ in our definition is a unit translation of the statistic $\widehat{R} = 1+ B/W$,
as defined in \cite{Gelman1996}.

We run a MH chain with local proposal specified by a uniform cube of side $6\sigma_1$ centered 
at the current state.  Using this LMHS we try to draw samples from the following needle in the haystack, 
i.e., $h$ with the following parameters:
\begin{equation}\label{E:NeedleParam}
\mu_1 = (0,0,0)^{\prime}, \mu_2 = (1,1,1)^{\prime}, \sigma_1 = 1, \sigma_2 = 0.006.
\end{equation}
 To diagnose convergence of the LMHS we calculate $B/W$ for each component of $x$ and
assume that the chain's burn-in time (the time when 
the samples may be affected by the initial condition) has ended when $B/W \leq 0.05$ for all three components.
The post burn-in run length, i.e., the number of samples kept after the burn-in, is set to be 
$100$ times the burn-in time (typical run lengths ranged in $[10000,50000]$ for target $h$ specified 
by \eqref{E:NeedleParam}).

\begin{figure}
\centering   \makebox{
\input{mix2norm_metro1.tex}
}
\caption{The running mean for four MH chains, as well as $B$, $W$, and $B/W$ for $x_1$ 
as a function of run length.  The true mean for $x_1$ is at $0.5$. \label{F:MHDiagnFail}}
\end{figure}
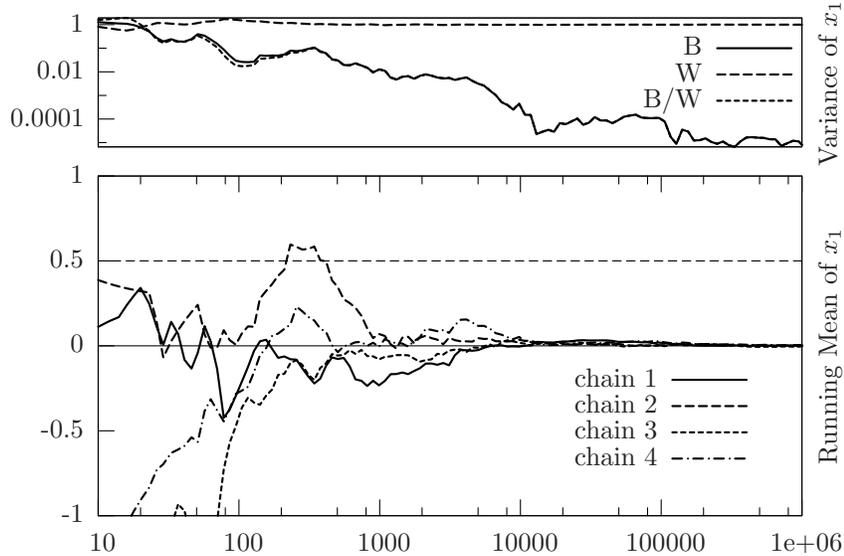

The above convergence diagnostics are more conservative than the standard recommendations
 \cite{Gelman1992,Gelman1996,Kass1998}.  Figure \ref{F:MHDiagnFail} shows the 
results (along the $x_1$ axis) of the above LMHS that relies on the $B/W$ statistic from four randomly 
initialized chains.  The running mean for each of the four chains has converged to the haystack mean 
of $(0,0,0)^{\prime}$ and completely missed the needle at $(1,1,1)^{\prime}$.  Thus, if we relied on our 
convergence diagnostic $B/W$, which appears to be consistently vanishing and thus suggestive of 
convergence to our target $h$, we would have entirely missed the needle.  Tuning the diagnostic 
parameters, including the number of chains, burn-in time, and run length, does not help diagnose true 
convergence for much sharper needles ($\sigma_2 < 10^{-5}$) that are naturally amenable to our MRS.

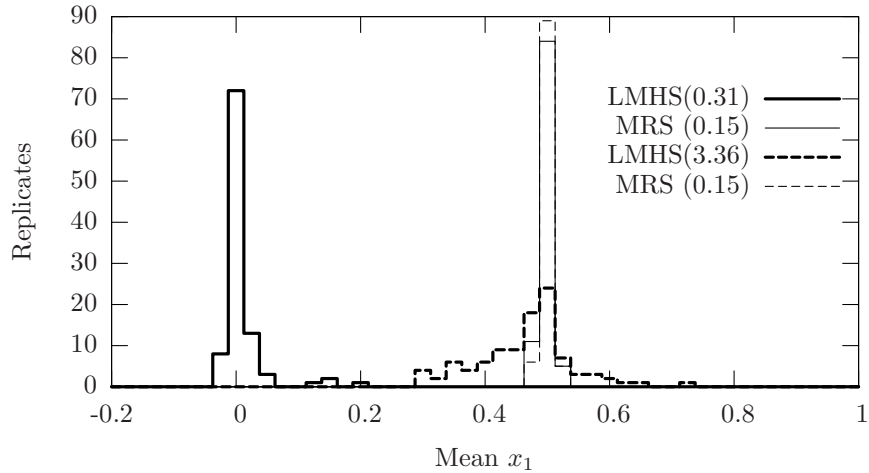
\begin{figure}
\centering   \makebox{
\input{mrsmetrocompare.tex}
}
\caption{Histograms of the mean $x_1$ from $100$ replicates of the LMHS 
and MRS.  The broken lines and solid lines represent targets with a blunt 
needle ($\sigma_2=0.10$) and a sharp needle ($\sigma_2=0.01$), respectively.  The CPU
time in seconds for each sampler is given in parenthesis. \label{F:LMRSMHS}}
\end{figure}

Next we compare the samples obtained from the $B/W$ diagnosed LMHS described above with 
$10,000$ samples from MRS induced by an integral-based adaptive partitioning of $\mathbf{\Theta}$ 
into $1,000$ boxes.  We compare the two samplers on two targets: (1) a blunt needle 
with $\sigma_2 = 0.10$ and (2) a sharp needle with $\sigma_2=0.01$.  The other parameters of the
two targets are the same as before \eqref{E:NeedleParam}.  The results are summarized in 
Figure \ref{F:LMRSMHS}.  The diagnostic $B/W$ works better in diagnosing convergence to the blunt 
target.  The bias is severe for the sharp needle in all $100$ replicates.  MRS clearly outperforms 
LMHS, both in terms of producing the true samples and in terms of CPU time 
(Figure \ref{F:LMRSMHS}).  Moreover, the sharpness of the needle only has a minor effect on the 
efficiency of MRS.  For example, for a much sharper needle with $\sigma_2 = 10^{-10}$, 
the MRS induced by an integral-based adaptive partitioning of $\mathbf{\Theta}$ into just $120$ cuboids, 
achieves an acceptance probability of 0.40.

\subsubsection{Rejection, Importance and Independent Metropolis-Hastings}

The same proposal density used in RS may be used as the proposal in importance 
sampler (IS) \cite{Kahn1953,Marshall1956} or as the proposal of the independent base chain in IMHS.  The 
latter two samplers are typically more efficient than RS,
although in some cases the efficiency of IMHS can be as low as half that
of RS \cite{Liu1995}.  The disadvantage 
of IMHS and IS (or RS) compared to MRS is in terms of diagnosing convergence and finding the right
proposal(s), respectively.  However, if one shares the proposal obtained through interval methods in MRS with IS and
IMHS, then we get their Moore versions which circumvent the disadvantages that arise from 
non-rigorously constructed proposals.  Indeed all three samples may be generated 
simultaneously from the same sequence of proposed values \cite{Cai1999}; each proposed value 
would be output with its importance weight, with some subset of the proposed values marked as 
IMHS-accepted, and with some further subset of those additionally marked as MRS-accepted and thereby 
constituting our collection of independent samples.  

\begin{figure}
\centering   \makebox[5pt]{
\input{RSMISIS.tex}
}
\caption{MSE of the three samplers, namely, MRS, IMHS and IS, as well as the acceptance probability 
of MRS and IMHS as a function of partition size (see text for description).
\label{F:3SamplersMSE}}
\end{figure}
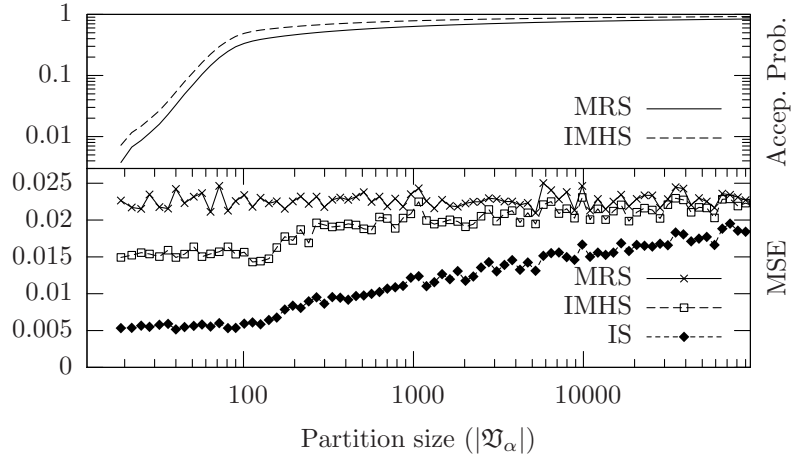

Figure \ref{F:3SamplersMSE} shows the mean squared error MSE for the sampler trio  
as a function of the size of the partition that is invoking their common proposal.  
The sample trio is drawn from our target $h$ \eqref{E:Needle} with the sharp needle ($\sigma_2=0.01$).  
To obtain the MSE for each sampler with target $p$ and proposal $q$, we drew $x_i \sim q, \, i = 1, \ldots, N$ using MRS, 
where $N$ is the number of samples needed to obtain 100 Moore rejection samples.  For IS each of the $x_i$'s were 
assigned the importance sampling weight $w_i = p(x_i)/q(x_i)$ and the estimated mean 
$\widehat{\mu} = {\sum_{i=1}^N{(w_i x_i)}}/{\sum_{i=1}^N{w_i}}$.  
The MRS estimated mean is $\widehat{\mu} = {\sum_{i=1}^{100}{x_{r_i}}}/100$, where $x_{r_i}$ is the $i^{th}$ MRS sample.  
For IMHS the mean is estimated by
$\widehat{\mu} = {\sum_{i=r_1}^N{x_i}}/{(N-r_1+1)}$, where $r_1$ is the index of the first MRS sample;
the early samples $x_i, \, i<r_1$ are excluded as burn-in.  This mean estimation was repeated $500$ times to 
obtain $\widehat{\mu}_j, \, j=1, \ldots, 500$ for each sampler.
Finally, the MSE was computed with 
the known mean $\mu=(0.5,0.5,0.5)$ under the Euclidean norm $|| \widehat{\mu}_j - \mu ||$ 
as $\sum_j{||\widehat{\mu}_j - \mu||^2/500}$.

The Figure \ref{F:3SamplersMSE} compares the three samplers and shows a typical 
pattern: at low acceptance probability, IS has lowest MSE, and MRS the highest, while at 
high acceptance probability all three samplers approach the same MSE. 
The lower MSE of IS is due to the large number of (MRS-discarded) samples being 
appropriately weighted.  Observe that such an auto-validating Moore importance 
sampler can be efficient and rigorous in estimating some expectation $E_p f(x)$ of interest.  
As the acceptance probability of MRS increases with refinement of the domain and the 
number of samples from each sampler approaches equality, the MSE of all three samplers converge as 
expected.  For some target shapes, e.g. the witches hat \eqref{E:WH}, we have observed the MSE
of IMHS to be greater than that of MRS, but by less than a factor of $2$,
in agreement with \cite{Liu1995} (results not shown).

\subsection{Multivariate Rosenbrock and Witch's Hat}\label{S:MRWH}

Next we examine the effect of dimensionality on efficiency of MRS through the 
challenging Rosenbrock function from the optimization literature.  We make it a Rosenbrock 
density ($r_D$) in $D$ dimensions over some compact $\mathbf{\Theta} \in \mathbb{IR}^D$ by 
appropriately normalizing the Rosenbrock shape $r^*_D$ \eqref{E:RosenDen}. 
\begin{equation}\label{E:RosenDen}
r^*_D (X) = \exp \{ - \sum_{i=2}^D (100 (X_i - X_{i-1}^2)^2 + (1-X_{i-1})^2 ) \}
\end{equation}

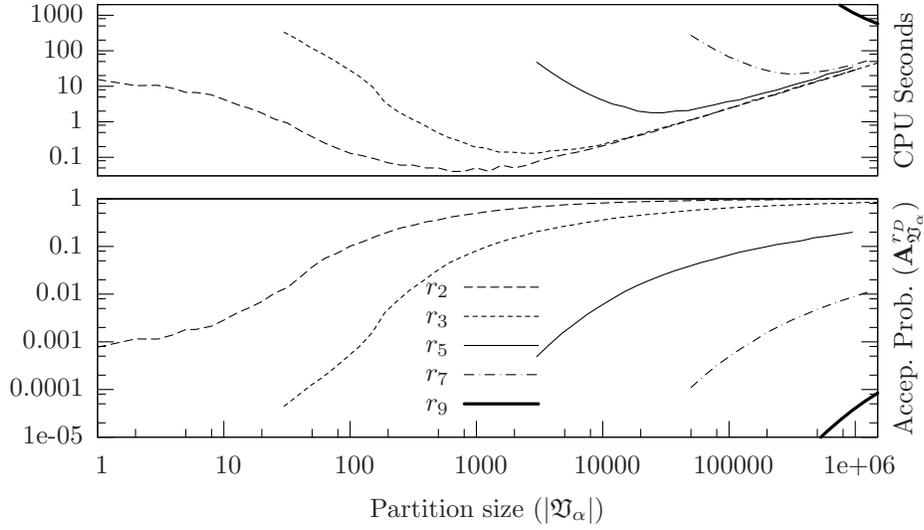
\begin{figure}
\centering   \makebox[5pt]{
\input{Rosen_compare_Temp.tex}
}
\caption{Acceptance probability ($\mathbf{A}^{r_D}_{\mathfrak{V}_{\alpha}}$) and CPU time  
to generate $10^4$ samples, as a function of partition size ($|\mathfrak{V}_{\alpha}|$), for Rosenbrock 
targets $r_D$ over $\mathbf{\Theta} = [-10,10]^{\otimes D}$, where $D$ is the dimension.
\label{F:RosenComp}}
\end{figure}

Figure \ref{F:RosenComp} summarizes the efficiency for various Rosenbrock densities.  For 
the more demanding nine dimensional Rosenbrock target $r_9$, we were able to draw $10,000$ 
samples in about $650$ CPU seconds at an acceptance probability of $1/10,000$.  The acceptance probability can 
be improved and/or $D$ can be increased naively if we allowed the partition size to be greater than 
a million.  Thus, the extent of RAM (random access memory) at our disposal ultimately determines 
the complexity and dimensionality of the target that can be rigorously sampled with MRS.  However,
the manner in which the natural interval extension is constructed will greatly affect the sampler's
efficiency as discussed later.  The acceptance probability for the relatively less complicated 
multivariate exponential mixture density truncated over 
$\mathbf{\Theta} = [-100,100]^{\otimes 10}$ is higher at $1/1000$ compared to that for the Rosenbrock
target $r_9$ even when there were $10$ modes inside a $10$-dimensional $\mathbf{\Theta}$ (results not shown).  
Thus, the complexity of the target greatly affects efficiency.

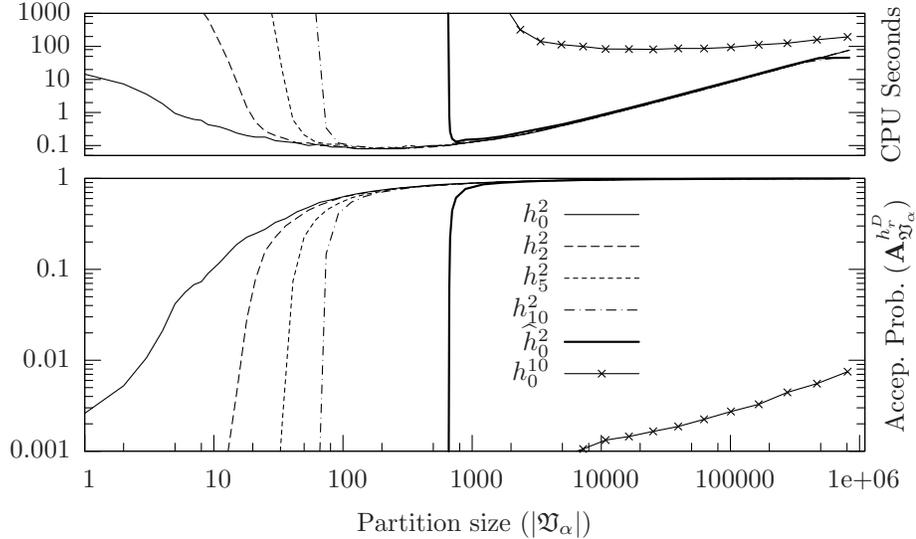
\begin{figure}
\input{Wh1_compare_Dim_New.tex}
\caption{Acceptance probability and CPU time  
to generate $10^4$ samples, versus partition size for witch's hat targets $w_r^D$, where $D$ is the dimension of the domain 
and $R = 10^{-r}$ is the hat's radius \eqref{E:WH}.  The hats of all targets were centered 
at the two vector $(2,\ldots,2)$.  The domain $\mathbf{\Theta}$ for  
$\widehat{h}_0^2$ was $[-10^{100}, +10^{100}]$, but all other targets had 
$\mathbf{\Theta} = [-10,10]^{\otimes D}$. \label{F:WH}}
\end{figure}

Finally, we arrive at the infamous witch's hat density which is considered to be a pathological 
target for most samplers \cite{Kass1998}.  The density is often thought of in two dimensions as an $m:(1-m)$ mixture of 
a cone with center $C$ and basal radius $R$ and a uniform distribution on a rectangle.  It can be 
easily generalized to $D$ dimensions as follows:
\begin{alignat}{2}\label{E:WH}
w_r^D(X) = & m \, \mathbf{1}_{ \{|| X-C || \leq R \} } \, \left( 1 - \frac{|| X-C ||}{R} \right) H + (1-m) \, \frac{1}{V}, \, \text{where} \notag \\
H =  & \frac{\Gamma(D/2) D(D+1)}{2 \pi^{D/2} R^D}, \quad V =  \prod_{i=1}^{D} {d(\mathbf{\Theta_i})}, \quad R = 10^{-r}. 
\end{alignat}
Our formulation of the witch's hat is even more challenging than the differentiable formulation suggested 
in \cite{Kass1998}, as the gradient is $0$ over the entire brim.  MRS is amenable to any target 
with a well-defined interval extension over the domain including $w_r^D$.  Mixtures of several sharply-peaked 
bivariate normals with  a uniform distribution, a further generalization of the other formulation \cite{Kass1998}, 
pose no sampling problems to MRS.  Figure \ref{F:WH} shows that one can efficiently sample 
from witch's hat  targets by rigorously constructing 
envelopes through the natural interval extension of \eqref{E:WH}.  We can even sample from the
hat of an eleven dimensional witch ($w_0^{10}$).  We can also make the
brim of the hat as large as $[-10^{100},+10^{100}]^{\otimes 2}$ without much trouble ($\widehat{w}_0^2$).  Note that 
decreasing the radius has a similar effect as widening the brim, in terms of lowering the acceptance
probability as a function of partition size.  Note that we are able to sample rigorously from a range of multivariate 
witch's hat targets with reasonable partition sizes and CPU seconds.

\subsection{Likelihood of Jukes-Cantor Triplets and Quartets}\label{S:JCTQ}

\begin{figure}
\makebox{\centerline{\includegraphics[width=5.50in]{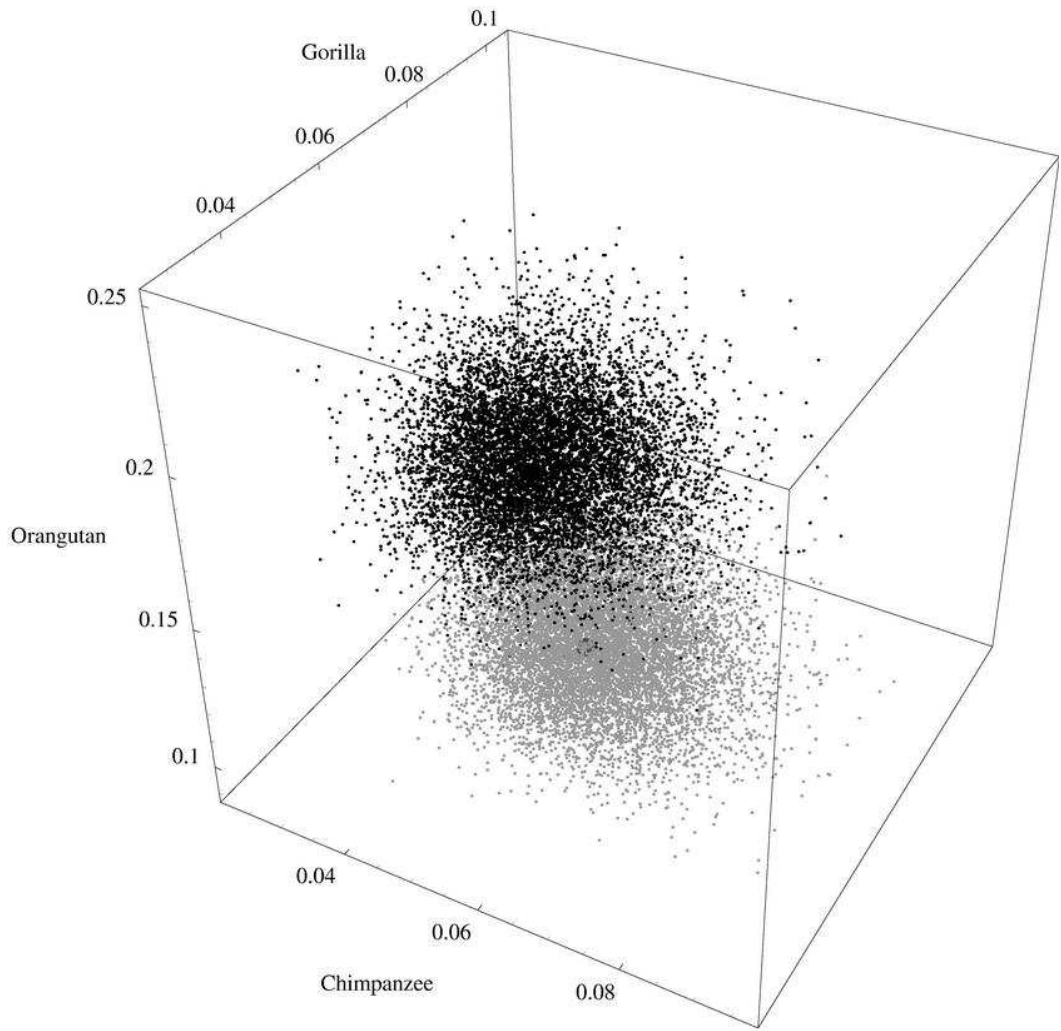}}} 
\caption{$10,000$ Moore rejection samples (gray dots) from the posterior distribution over the three branch 
lengths of the unrooted phylogenetic tree space of Chimpanzee, Gorilla and Orangutan based on their 
mitochondrial DNA.  Marginal triplets (dark dots) of $10,000$ samples from the quartet tree space 
of Chimp, Gorilla, Orangutan and Gibbon.  \label{JC3MRS}}
\end{figure}

Inferring the ancestral relationship among a set of species based on their DNA sequences is a
basic problem in phylogenetics \cite{Semple2003,Felsenstein03}.  One can obtain the likelihood of a particular
phylogenetic tree that relates the species of interest by superimposing a simple Markov model of
DNA substitution due to Jukes and Cantor \cite{Jukes1969} on that tree.  The length of an edge
(branch length) connecting two nodes (species) in the tree represents the amount of evolutionary time 
(divergence) between the two species.  The likelihood function over trees obtained 
through a post-order traversal (e.g.~\cite{Felsenstein1981}) has a natural interval extension over boxes 
of trees \cite{SainudiinYoshida2005}.  This allows us to draw samples from the posterior distribution over
some compact box in the tree space using our MRS.  Using the data from the mitochondria of 
Chimpanzee, Gorilla, and Orangutan \cite{Brown1982} that can be summarized by 29 distinct site 
patterns \cite{Sainudiin2004}, we 
obtain the posterior distribution by normalizing the likelihood with a uniform prior over the biologically meaningful 
compact domain $\mathbf{\Theta} = [10^{-10}, 10]^{\otimes 3}$.  $10,000$ independent samples were drawn 
in $942$ CPU seconds from the posterior distribution over Jukes-Cantor triplets, i.e. unrooted trees with three 
edges corresponding to the three primates emanating from their common ancestor.  Figure \ref{JC3MRS} shows 
these samples (gray points) scattered about the verified global MLE of the triplet \cite{SainudiinYoshida2005}.  

We were able to draw samples from Jukes-Cantor quartets (unrooted trees for four taxa) by adding the homologous 
sequence of the Gibbon which resulted in 61 distinct site patterns \cite{Sainudiin2004}.  This is a more 
challenging problem because there are 3 distinct tree topologies for an unrooted quartet tree and 
each of these has five edges.  Thus, the domain of quartets is a piecewise 
Euclidean space that arises from a fusion of 3 distinct five dimensional orthants.  Since the post-order
traversals specifying the likelihood function are topology-specific, we extended the likelihood over 
a compact box of quartets in a topology-specific manner.  The computational
time was about a day and a half to draw $10,000$ samples from the quartet target due to low 
acceptance probability of the naive likelihood function based on distinct site patterns.  All the samples had 
the same topology which grouped Chimp and Gorilla together, i.e. ((Chimp, Gorilla), (Orangutan, Gibbon)).  The 
samples were again scattered about the verified global MLE of the quartet \cite{Sainudiin2004}.  The marginal 
triplet trees (dark dots) within the $10,000$ sampled quartets are also plotted in Figure \ref{JC3MRS}.  Observe 
the influence of an additional taxon on the triplet estimates.
This quartet likelihood function has an elaborate DAG (Definition \ref{D:DAG}) with numerous operations.  
When the data got compressed into sufficient statistic through algebraic statistical methods 
\cite{PachterSturmfels2005}, the efficiency increased tremendously (for e.g.~triplet efficiency increases by a 
factor of $3.7$).  This is due to the 
number of leaf nodes in the target DAG, which encode the distinct site patterns of the observed data 
into the likelihood function, getting reduced from $29$ to $5$ for the triplet target and from 
$61$ to $15$ for the quartet target \cite{Casanellas2005}.  Poor sampler efficiency makes it impractical to 
sample from trees with
five or more leaves.  However, one could use such triplets and quartets drawn from the 
posterior distribution to stochastically amalgamate and produce estimates of larger trees via 
fast amalgamating algorithms \cite{Strimmer1996,Levy2005}.  A collection 
of large trees obtained through such stochastic amalgamations would account for the effect of finite
sample sizes (sequence length) as well as the sensitivity of the amalgamating algorithm itself 
to variation in the input vector of small tree estimates.  It would be interesting to investigate if such
stochastic amalgamations can help improve mixing of MCMC algorithms on large tree spaces \cite{Mossel2005}.

\section{Conclusion} \label{S:C}

Interval methods provide a rigorous, efficient and fairly general way of constructing envelope 
functions for use in rejection sampling from target densities with a well-defined interval extension.  
In particular the method allows the envelope to be drawn from a large, flexible family of functions 
(simple functions over a family of adaptively refined partitions), and to be constructed in a manner that
rigorously maintains the envelope property as the envelope function is adaptively refined.  Refining the 
partition decreases the rejection probability at a rate that is no slower than linear with the mesh.  The 
corresponding proposal density is easily constructed in $\mathcal{O}$(partition size) time into a data 
structure that allows samples from it to be drawn in constant time.  When one substitutes conventional 
floating-point arithmetic for real arithmetic in a computer and uses discrete lattices to construct the envelope 
and/or proposal, it is generally not possible to guarantee the envelope property and thereby ensure that 
samples are drawn from the desired target density,
except in special cases.  For example, the adaptive rejection sampler (ARS) \cite{Gilks1992,GilksWild1992} 
is efficient at drawing independent samples only from one dimensional log-concave targets.  In  
ARS as well as a subsequent generalization of it through a Metropolis sampling step \cite{Gilks1995} to one 
dimensional non-log-concave targets, one can draw samples from higher 
dimensional targets by Gibbs sampling \cite{Geman1984,Gelfand1990} one dimension 
at a time.  On one hand, Gibbs sampling, being
a special case of Metropolis-Hastings sampling \cite{Chib1995}, is at the mercy of heuristic convergence
diagnostics.  On the other, proposals constructed for non-log-concave conditional densities from finitely 
many points cannot guarantee that the density has not soared between the sampled points.
However, the construction of the Moore rejection sampler through interval methods, that enclose the target shape 
over the entire real continuum in any box of the domain with machine-representable bounds, in 
a manner that rigorously accounts for all sources of numerical errors (see \cite{Kulisch2001,Hammer1995} 
for a discussion on error control), naturally guarantees that the Moore rejection samples are independent 
draws from the desired target.  Moreover, the target is allowed to be multivariate and/or non-log-concave 
with possibly `pathological' behavior, as long as it has a well-defined interval extension.  

Unfortunately, the efficiency of MRS is not immune to the curse of dimensionality and target DAG complexity.  
When the DAG
for the likelihood gets large, its natural interval extension can have terrible over-enclosures of the true
range, which in turn forces the adaptive refinement of the domain to be extremely 
fine for efficient envelope construction.  Thus, a naive application of interval methods to targets 
with large DAGs can be terribly inefficient.  In such cases, sampler efficiency rather than rigor is the 
issue.  Thus, one will not obtain samples in a reasonable time rather than produce samples from some 
unknown and undesired target.  There are several ways in which efficiency can be improved for such 
cases.  First, the particular structure of the target DAG should be exploited to avoid any redundant computations.
For example, algebraic statistical methods can be used to find sufficient statistics to dissolve symmetries 
in the DAG as done in Section \ref{S:JCTQ}.  
Second, we can further improve efficiency by limiting ourselves to differentiable targets in $C^n$.
Tighter enclosures of the range $p^*(\Theta^{(i)})$ with $P^*(\Theta^{(i)})$ can come from the 
enclosures of Taylor expansions of $p^*$ around the midpoint $m(\Theta^{(i)})$ 
through interval-extended automatic differentiation \cite{Rall1981,Berz1991,Kulisch2001} that can
then yield tighter estimates of the integral enclosures \cite{Tucker2004}. 
Third, we can employ pre-processing to improve efficiency.  For example, we can pre-enclose the range of a possibly 
rescaled $p^*$ over a partition of the domain and then obtain the 
enclosure of $P^*$ over some arbitrary $\Theta$ through a combination of hash access and hull 
operations on the pre-enclosures.  Such a pre-enclosing technique reduces not only the overestimation
of target shapes with large DAGs but also the computational cost incurred while 
performing interval operations with processors that are optimized for floating-point arithmetic.  Fourth, efficiency at the 
possible cost of rigor can also be gained  (up to $30\%$ ) by foregoing directed rounding during envelope 
construction.  

In this paper we focused on the interval extension of the simplest sampler, namely the rejection sampler.  
We also exploited the direct connections between rejection sampler, importance sampler 
and independent Metropolis-Hastings sampler to produce sample trios from the interval extensions of all three
samplers.  It would be interesting to compare other Monte Carlo methods
to their natural interval extensions.  For example, even Metropolis-coupled MCMC \cite{Geyer1991} which 
was designed to accelerate convergence for complicated targets is known to converge exponentially slowly in some
cases \cite{Bhatnagar2004}.  Preliminary analysis suggests that a non-rigorous interval extension of the 
local Metropolis-Hastings sampler (relying on heuristic convergence diagnostics) 
may have a higher probability of converging to the target when compared to its 
floating-point cousin.  Such hybrid samplers that rely on both interval and local methods may efficiently 
produce fairly reliable samples from challenging higher dimensional targets.

\section{Acknowledgments}
This was supported by a joint NSF/NIGMS grant DMS-02-01037 to Durrett, Aquadro, and Nielsen.  R.S.~is
a Research Fellow of the Royal Commission for the Exhibition of 1851.
Many thanks to Rob Strawderman and Warwick Tucker for constructive comments.

\section{Appendix A}\label{S:AppA}
\subsubsection*{Arithmetic on intervals in $\mathbb{IR} \triangleq \{ [x,y] : x \leq y, x,y \in \mathbb{R} \}$}

\begin{definition}[Interval arithmetic]\label{Df:intarith}
If the binary operator $\star$ is one of the elementary arithmetic 
operations $\{ +,-,\cdot,/\}$, then we define an arithmetic on operands in 
$\mathbb{IR}$ by
\[
X \star Y \triangleq \{x \star y : x \in X , y \in Y\}
\]
with the exception that $X/Y$ is undefined if $0 \in Y$.
\end{definition}
\begin{theorem}\label{Th:IRarith}
Arithmetic on the pair $X,Y \in \mathbb{IR}$ is given by:
\[
\begin{array}{lcl}
X + Y  & = & [\underline{x} + \underline{y}, \overline{x} + \overline{y}] \\
X - Y & = & [\underline{x} - \overline{y}, \overline{x} - \underline{y}] \\
X \cdot Y & = & [\min \{\underline{x} \underline{y}, \underline{x} \overline{y},
\overline{x} \underline{y}, \overline{x} \overline{y}\},
\max \{\underline{x} \underline{y}, \underline{x} \overline{y},
\overline{x} \underline{y}, \overline{x} \overline{y}\}], \notag \\
X / Y & = & X \cdot [1/\overline{y}, 1/\underline{y}], \text{ provided, } 0 \notin Y. \notag
\end{array}
\]
\end{theorem}
{\bf Proof} (cf.~\cite{Hammer1995,Tucker2004}):
Since any real arithmetic operation $x \star y$, where $\star \in \{ +,-,\cdot,/\}$ and
$x,y \in \mathbb{R}$, is a continuous function
$x \star y\triangleq\star(x,y): \mathbb{R} \times \mathbb{R} \rightarrow \mathbb{R}$,
except when $y=0$ under $/$ operation.  Since $X$ and $Y$ are simply connected compact 
intervals, so is their product $X \times Y$.  On such a domain $X \times Y$, the continuity 
of $\star(x,y)$ (except when $\star = /$ and $0 \in Y$) ensures the attainment of a minimum, 
a maximum and all intermediate values.  Therefore, with the exception of the case when
$\star = /$ and $0 \in Y$, the range $X \star Y$ has an interval form
$[\min{(x \star y)}, \max{(x \star y)}]$, where the $\min$ and $\max$ are taken over all
pairs $(x,y) \in X \times Y$.  Fortunately, we do not have to evaluate $x \star y$ over 
every $(x, y) \in X \times Y$ to find the global $\min$ and global $\max$ of $\star(x,y)$ 
over $X \times Y$, because the monotonicity of the $\star(x,y^*)$ in terms of $x \in X$ for 
any fixed $y^* \in Y$ implies that the extremal values are attained on the boundary of 
$X \times Y$, i.e., the set $\{\underline{x}, \underline{y}, \overline{x}$, and $\overline{y} \}$.  
Thus the theorem can be verified by examining the finitely many boundary cases. $\square$

An extremely useful property of interval arithmetic that is a direct consequence of 
Definition \ref{Df:intarith} is summarized by the following theorem.  
\begin{theorem}[Fundamental property of interval arithmetic]\label{Th:FundPropIntArith}
If $X \subseteq X^{\prime}$ and $Y \subseteq Y^{\prime}$ and $\star \in \{+,-,\cdot, / \}$, then
\[
X \star Y \subseteq X^{\prime} \star Y^{\prime}, 
\]
where we require that $0 \notin Y^{\prime}$ when $\star = /$.
\end{theorem}
{\bf Proof}:
\[
X \star Y = \{ x \star y: x \in X, y \in Y \} \subseteq \{ x \star y: x \in X^{\prime}, y \in Y^{\prime} \} = 
X^{\prime} \star Y^{\prime}. \square
\]
Note that an immediate implication of Theorem \ref{Th:FundPropIntArith} is that when $X = x$ and $Y = y$ are thin
intervals (real numbers $x$ and $y$), then $X^{\prime} \star Y^{\prime}$ will contain the result of the real arithmetic 
operation $x \star y$.
\begin{definition}[Range]
Consider a real-valued function $f: D \rightarrow \mathbb{R}$ where the domain $D \subseteq \mathbb{R}^n$.  
The range of $f$ over any $E \subseteq D$ is represented by $Rng(f;E)$ and defined to be the set
\[
Rng(f;E) \triangleq \{ f(x) : x \in E \}
\]
However, when the range of $f$ over any $X \in \mathbb{IR}^n$ such that $X \subseteq D$ is of interest, 
we will use the short-hand $f(X)$ for $Rng(f;X)$.
\end{definition}

\begin{definition}[Interval extension of subsets of $\mathbb{R}^n$] \label{D:Boxes}
For any Euclidean subset $\mathbf{\Theta} \subseteq \mathbb{R}^n$ let us denote its interval extension 
by $\mathbb{I}\mathbf{\Theta}$ and define it to be the set
\[
\mathbb{I}\mathbf{\Theta} \triangleq \{ X \in \mathbb{IR}^n : \underline{x}, \overline{x} \in \mathbf{\Theta} \}
\]
We refer the the $k$th interval of interval vector or box $X \in \mathbb{IR}^n$ by $X_k$.
\end{definition}

\begin{definition}[Inclusion isotony]
An box-valued map $F: D \rightarrow \mathbb{IR}^m$, where $D \in \mathbb{IR}^n$, is inclusion isotonic if it 
satisfies the property 
\[
\forall \, X \subseteq Y \subseteq D \implies F (X) \subseteq F(Y).
\]
\end{definition}

\begin{definition}[The natural interval extension]\label{D:NIE}
Consider a real-valued function $f: D \rightarrow \mathbb{R}$ given by a formula, where the 
domain $D \in \mathbb{IR}^n$.  If real constants, variables, and operations in $f$ are replaced 
by their interval counterparts, then one obtains
\[
F(X): \mathbb{I}D \rightarrow \mathbb{IR}.
\]
$F$ is known as the natural interval extension of $f$.  This extension is well-defined if we do not run into
division by zero.
\end{definition}

\begin{theorem}[Inclusion isotony of rational functions]\label{Thm:Rational}
Consider the rational function $f(x) = p(x)/q(x)$, where $p$ and $q$ are polynomials.  Let $F$ be its 
natural interval extension such that $F(Y)$ is well-defined for some $Y \in \mathbb{IR}$ and let
$X, X^{\prime} \in \mathbb{IR}$.  Then we have
\[
\begin{array}{ccl}
(i) & \text{\em Inclusion isotony:} & \forall \, X \subseteq X^{\prime} \subseteq Y \implies F(X) \subseteq F(X^{\prime}) 
\, \text{\em , and } \\ 
(ii) & \text{\em Range enclosure:}   & \forall \, X \subseteq Y \implies Rng(f; X) = f(X) \subseteq F(X).
\end{array}
\]
\end{theorem}
{\bf Proof} (cf.~\cite{Tucker2004}):
Since $F(Y)$ is well-defined, we will not run into division by zero, and therefore (i) follows from the repeated 
invocation of Theorem \ref{Th:FundPropIntArith}.  We can prove (ii) by contradiction.  Suppose $Rng(f; X) \nsubseteq F(X)$.  
Then there exists $x \in X$, such that $f(x) \in Rng(f;X)$ but $f(x) \notin F(X)$.  This in turn implies that 
$f(x)=F([x,x]) \notin F(X)$, which contradicts (i).  Therefore, our supposition cannot be true and we have 
proved (ii) $Rng(f; X) \subseteq F(X)$. $\square$  

\begin{definition}[Standard functions]
Piece-wise monotone functions, including exponential, logarithm, rational power, absolute value, and trigonometric functions,
constitute the set of standard functions
\[
\mathfrak{S} = \{ \, a^x, {log}_b(x), x^{p/q}, |x|, \sin(x), \cos(x), \tan(x), \sinh(x), \ldots, \arcsin(x), \ldots \, \}.
\]
\end{definition}
Such functions have well-defined interval extensions that satisfy inclusion isotony and {\em exact range enclosure}, i.e., 
$Rng(f;X) = f(X) = F(X)$.  Consider the following definitions for the interval extensions for some monotone functions 
in $\mathfrak{S}$ with $X \in \mathbb{IR}$,
\[
\begin{array}{lclr}
\exp(X)		&=& [\exp(\underline{x}), \exp(\overline{x})] 	& \\
\arctan(X)	&=& [\arctan(\underline{x}), \arctan(\overline{x})] 	& \\
\sqrt{(X)}	&=& [\sqrt{(\underline{x})}, \sqrt{(\overline{x})}] & \text{ if } 0 \leq \underline{x} \\
\log(X)		&=& [\log(\underline{x}), \log(\overline{x})]   & \text{ if } 0 < \underline{x} \\
\end{array}
\] 
and a piece-wise monotone function in $\mathfrak{S}$ with $\mathbb{Z}^+$ and $\mathbb{Z}^-$ representing the 
set of positive and negative integers, respectively.
\begin{equation}
X^n = 
\begin{cases}
[\underline{x}^n,\overline{x}^n] & \text{:   if } n \in \mathbb{Z}^+ \text{ is odd}, \notag \\
[{\langle X \rangle}^n,{|X|}^n] & \text{:   if } n \in \mathbb{Z}^+ \text{ is even}, \notag \\
[1,1] & \text{:   if } n = 0, \notag \\
[1/ \overline{x},1/ \underline{x}]^{-n} & \text{:   if } n \in \mathbb{Z}^- ; 0 \notin X \notag 
\end{cases}
\end{equation}

\begin{definition}[Elementary functions]\label{D:ElemFunc}
A real-valued function that can be expressed as a finite combination of constants, variables, arithmetic operations, 
standard functions and compositions is called an elementary function.  The set of all such elementary functions is
referred to as $\mathfrak{E}$.
\end{definition}

\begin{definition}[Directed acyclic graph (DAG) of a function]\label{D:DAG}
One can think of the process by which an elementary function $f$ is computed as the result of
a sequence of recursive operations with the subexpressions $f_i$ of $f$ where, $i=1,\ldots,n < \infty$.  This 
involves the evaluation of the subexpression $f_i$ at node $i$ with operands $s_{i_i},s_{i_2}$ from 
the sub-terminal nodes of $i$ given by the directed acyclic graph (DAG) for $f$
\begin{equation}\label{E:odotf}
s_i = \odot f_i \triangleq 
\begin{cases}
f_i(s_{i_1},s_{i_2}) & \text{:  if node $i$ has 2 sub-terminal nodes $s_{i_1},s_{i_2}$}  \\ 
f_i(s_{i_1}) & \text{:  if node $i$ has 1 sub-terminal node $s_{i_1}$}   \\ 
I(s_{i}) & \text{:  if node $i$ is a leaf or terminal node, } $I(x) = x$.   
\end{cases}
\end{equation}
The leaf or terminal node of the DAG is a constant or 
a variable and thus the $f_i$ for a leaf $i$ is set equal to the respective constant or variable.  The 
recursion starts at the leaves and terminates at the root of the DAG.  The DAG for an 
elementary $f$ with $n$ sub-expressions $f_1,f_2,\ldots,f_n$ is :
{\large 
\begin{equation}\label{E:DAGRec}
\{ \odot f_i \}_{i=1}^n \quad \rightarrowtail \quad \odot f_n = f(x), 
\end{equation}
}
where each $\odot f_i$ is computed according to \eqref{E:odotf}.
\end{definition}

\begin{figure}
   \makebox{\centerline{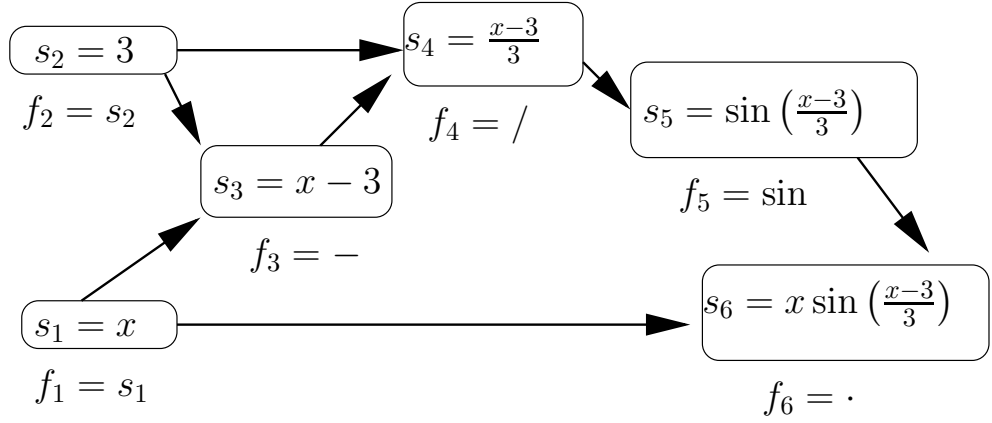}}
   \caption{Recursive evaluation of the sub-expressions $f_1,\ldots,f_6$ on the DAG
of the elementary function $f(x) = \odot f_6 =  x \cdot \sin ((x-3)/3)$}\label{Fi:DAG}
\end{figure}

For example the elementary function $x \cdot \sin ((x-3)/3)$
can be obtained from the terminus $\odot f_6$ of the recursion $\{ \odot f_i \}_{i=1}^6$ on the DAG for $f$ as shown 
in Figure \ref{Fi:DAG}.  
It would be convenient if guaranteed enclosures of the range $f(X)$ of an elementary $f$ can be 
obtained by its natural interval extension $F(X)$.  We show that inclusion
isotony does indeed hold for $F$, i.e.~if $X \subseteq Y$, then $F(X) \subseteq F(Y)$, and in
particular, the {\em inclusion property} that $x \in X \implies f(x) \in F(X)$ does hold.  

\begin{theorem}[The fundamental theorem of interval analysis]\label{3.1.11}
Consider any elementary function $f \in \mathfrak{E}$.  Let $F: Y \rightarrow \mathbb{IR}$ be its 
natural interval extension such that $F(Y)$ is well-defined for some $Y \in \mathbb{IR}$ and let 
$X, X^{\prime} \in \mathbb{IR}$.  Then we have
\[
\begin{array}{ccl}
(i) & \text{\em Inclusion isotony:} & \forall \, X \subseteq X^{\prime} \subseteq Y \implies F(X) \subseteq F(X^{\prime}) 
\, \text{\em , and } \\ 
(ii) & \text{\em Range enclosure:}   & \forall \, X \subseteq Y \implies Rng(f; X) = f(X) \subseteq F(X).
\end{array}
\]
\end{theorem}
{\bf Proof} (cf.~\cite{Tucker2004}): 
Any elementary function $f \in \mathfrak{E}$ is defined by the recursion \ref{E:DAGRec} on its 
sub-expressions $f_i$ where $i \in \{ 1,\ldots,n \}$ according to its DAG.
If $f(x) = p(x)/q(x)$ is a rational function, then the theorem already holds by Theorem \ref{Thm:Rational}, 
and if $f \in \mathfrak{S}$ then the theorem holds because the range 
enclosure is exact for standard functions.  Thus it suffices to show that if the theorem 
holds for $f_1, f_2 \in \mathfrak{E}$, then the theorem also holds for $f_1 \star f_2$, where 
$\star \in  \{ +,-,/,\cdot,\circ \}$.  By $\circ$ we mean the composition operator.  Since 
the proof is analogous for all five operators, we only focus on the $\circ$ operator.
Since $F$ is well-defined on its domain $Y$, neither the real-valued $f$ nor any of its 
sub-expressions $f_i$ have singularities in its respective domain $Y_i$ induced by $Y$.  
In particular $f_2$ is continuous on any $X_2$ and $X^{\prime}_2$ such that 
$X_2 \subseteq X^{\prime}_2 \subseteq Y_2$ implying the compactness of $F_2(X_2) \triangleq W_2$ and 
$F_2(X_2^{\prime}) \triangleq W_2^{\prime}$, respectively.  By our assumption that $F_1$ and $F_2$ 
are inclusion isotonic we have that $W_2 \subseteq W_2^{\prime}$ and also that 
\[
F_1 \circ F_2 (X_2) = F_1(F_2(X_2)) = F_1(W_2) \subseteq F_1(W_2^{\prime}) = F_1(F_2(X_2^{\prime})) = F_1 \circ F_2 (X_2)
\] 
The range enclosure is a consequence of inclusion isotony by an argument identical to that given in
the proof for Theorem \ref{Thm:Rational}. $\square$ 

The fundamental implication of the above theorem is that it allows us to enclose the range of
any elementary function and thereby produces an upper bound for the global maximum and a
lower bound for the global minimum over any compact subset of the domain upon which the function 
is well-defined.  We will see in the sequel that this is the work-horse of randomized 
enclosure algorithms that efficiently produce samples even from highly multi-modal target distributions.

\begin{figure}
\centering   \makebox{\includegraphics[width=5.0in]{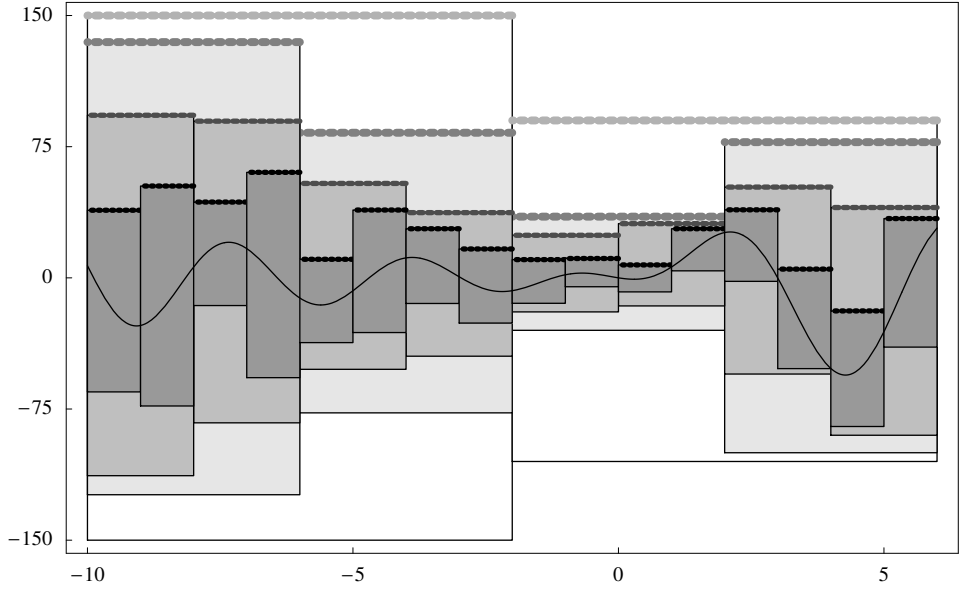}}  
   \caption{Range enclosure of the interval extension of 
   $-\sum_{k=1}^5{ k \, x \, \sin{(\frac{k(x-3)}{3})}}$ linearly tightens with the mesh.}\label{Fi:refine}
\end{figure}

Unlike the natural interval extension of an $f \in \mathfrak{S}$ that produces exact 
range enclosures, the natural interval extension $F(X)$ of an $f \in \mathfrak{E}$ often 
overestimates the range $f(X)$, but can be shown under mild conditions
to linearly approach the range as the maximal diameter of the box $X$ goes to zero, i.e.,
$\mathfrak{h}(F(X),f(X)) \leq \alpha \cdot d_{\infty}(X)$ for some $\alpha \geq 0$.  This
implies that a partition of $X$ into smaller boxes $\{X^{(1)}, \cdots , X^{(m)}\}$ gives better 
enclosures of $f(X)$ through the union
$\bigcup_{i=1}^m F(X^{(i)})$ as illustrated in Figure \ref{Fi:refine}.  
Next we make the above statements precise.

\begin{definition}\label{D:LipElemFunc}
A function $f:D \rightarrow \mathbb{R}$ is Lipschitz if there exists
a Lipschitz constant $K$ such that, for all $x,y \in D$, we have $|f(x)-f(y)| \leq K |x-y|$.
We define $\mathfrak{E_L}$ to be the set of elementary functions whose sub-expressions $f_i$, $i=1,\ldots,n$ 
at the nodes of the corresponding DAGs are all Lipschitz.  
\end{definition}

\begin{theorem}[Range enclosure tightens linearly with mesh]\label{3.1.15}
Consider a function $f:D \rightarrow \mathbb{R}$ with $f \in \mathfrak{E_L}$.  Let $F$ be an inclusion isotonic interval extension
of $f$ such that $F(X)$ is well-defined for some $X \in \mathbb{IR}, X \subseteq I$.  Then there exists a 
positive real number $K$, depending on $F$ and $X$, such that if $X = \cup_{i=1}^k X^{(i)}$, then
\[
Rng(f;X) \subseteq \bigcup_{i=1}^k F(X^{(i)}) \subseteq F(X)
\]
and
\[
r \left( \bigcup_{i=1}^k F(X^{(i)}) \right) \leq r(Rng(f;X)) + K \max_{i=1,\ldots,k} {r(X^{(i)})}
\]
\end{theorem}
{\bf Proof }:  The proof is given by an induction on the DAG for $f$ similar to the proof of 
Theorem \ref{3.1.11} (See \cite{Tucker2004}).

\section{Appendix B}\label{S:AppB}
Here we will study the Moore rejection sampler (MRS) carefully.  Lemma \ref{lemma1} 
shows that MRS indeed produces independent samples from the desired target and 
Lemma \ref{lemma2} describes the asymptotics of the acceptance probability as the partition
of the domain is refined.

\begin{lemma}\label{lemma1}
Suppose that the target shape $p^*$ has a well-defined natural interval extension $P^*$.
If $U$ is generated according to the steps in part \ref{Alg:RS} of the rejection sampling algorithm, 
and if the proposal density $q^{\mathfrak{T}}(\theta)$  and the envelope function $f_{q^{\mathfrak{T}}}(\theta)$ 
are given by \eqref{E:qRS} and \eqref{E:fqRS}, respectively, then $U$ is distributed according to the target $p$.
\end{lemma}
{\bf Proof}:
From  \eqref{E:qRS} and \eqref{E:fqRS} observe that $f_{q^{\mathfrak{T}}}(t) = q^{\mathfrak{T}}(t) N_{q^{\mathfrak{T}}}$.
Let us define the following two subsets of $\mathbb{R}^2$,
\[
\mathcal{B}_q = \{(t,h) : 0 \leq h \leq f_{q^{\mathfrak{T}}}(t) \}, \text{   and   } 
\mathcal{B}_p = \{(t,h) : 0 \leq h \leq p^*(t) \}.
\]
First let us agree that steps \ref{RS1} and \ref{RS2} of part \ref{Alg:RS} of the rejection sampling algorithm
produce a pair $(T,H)$ that is uniformly distributed on $\mathcal{B}_q$.  We can see this by
letting $k(t,h)$ denote the joint density of $(T,H)$ and $k(h | t)$ denote the conditional
density of $H$ given $T=t$.  Then,
\[
k(t,h) = 
\begin{cases}
q^{\mathfrak{T}}(t) \, k(h | t) & \text{ if } (t,h) \in \mathcal{B}_q \\
0   & \text{ otherwise }.
\end{cases}
\]
Since we sample a uniform height $h$ for a given $t$ in Step \ref{RS2} of the algorithm
\[
k(h | t) = 
\begin{cases}
(f_{q^{\mathfrak{T}}}(t))^{-1} = (q^{\mathfrak{T}}(t)  N_{q^{\mathfrak{T}}})^{-1} 
& \text{ if } h \in [0, f_{q^{\mathfrak{T}}}(t)] \\
0  & \text{ otherwise}.
\end{cases}
\] 
Therefore,
\[
k(t,h) = 
\begin{cases}
q^{\mathfrak{T}}(t) \, k(h | t) = q^{\mathfrak{T}}(t) / (q^{\mathfrak{T}}(t) \, N_{q^{\mathfrak{T}}} )
= (N_{q^{\mathfrak{T}}})^{-1} & \text{ if } (t,h) \in \mathcal{B}_q \\
0   & \text{ otherwise }.
\end{cases}
\]
Thus we have shown that the joint density of $(T,H)$ is a uniformly distribution on $\mathcal{B}_q$.  The 
above relationship also makes geometric sense since the volume of $\mathcal{B}_q$ is exactly $N_{q^{\mathfrak{T}}}$.
Now, let $(T^*,H^*)$ be an accepted point, i.e., $(T^*,H^*) \in \mathcal{B}_p \subseteq \mathcal{B}_q$.  Then, the
uniform distribution of $(T,H)$ on $\mathcal{B}_q$ implies the uniform distribution of $(T^*,H^*)$ on $\mathcal{B}_p$.  
Since the volume of $\mathcal{B}_p$ is $N_p$, the p.d.f.~of $(T^*,H^*)$ is identically $1/N_p$ on $\mathcal{B}_p$ and 
$0$ elsewhere.  Hence, the marginal p.d.f. of $U=T^*$ is
\[
\begin{array}{lcl}
w(u) & = & \int_0^{p^*(u)} 1/N_p \, d h \\ 
     & = & 1/N_p \int_0^{p^*(u)} 1 \, d h \\
     & = & 1/N_p \int_0^{N_p p(u)} 1 \, d h, \quad \because \, p(u) = p^*(u)/N_p \\
     & = & p(u). \quad \square  
\end{array}
\]
  
\begin{lemma} \label{lemma2}
Let $\mathfrak{U}_W$ be the uniform partition of $\mathbf{\Theta} = [\underline{\theta},\overline{\theta}]$ 
into $W$ intervals each of diameter $w$
\[
\begin{array}{lcl}
w & = &  \frac{(\overline{\theta} - \underline{\theta})}{W} \\ 
\Theta_W^{(i)} & = & [ \ \underline{\theta} + (i-1) w, \ \underline{\theta} + i w \ ] \, , i = 1, \dots, W \\ 
\mathfrak{U}_W & = & \{ \Theta_W^{(i)} \, , i = 1, \dots, W \}. 
\end{array}
\]
and let $p^* \in \mathfrak{E_L}$, then
\[
\mathbf{A}^p_{\mathfrak{U}_W} = 1 - \mathcal{O} (1/W)
\]
\end{lemma}
{\bf Proof} \\
Then by means of Theorem \ref{3.1.15} 
\[
\begin{array}{lcl}
d(\Theta_W^{(i)}) = \mathcal{O} (1/W) & \implies & \mathfrak{h}( \ p^*(\Theta_W^{(i)}), P^*(\Theta_W^{(i)}) \ ) 
= \mathcal{O} (1/W) \\
& \implies & d(P^*(\Theta_W^{(i)}) )=  \mathcal{O} (1/W), \qquad \because p^* \in {\mathfrak{E_L}}
\end{array}
\]
Therefore
\[
  \sum_{i=1}^{|\mathfrak{U}_W|} \left( d(\Theta_W^{(i)}) \cdot P^*(\Theta_W^{(i)}) \right) 
= w \sum_{i=1}^{W} P^* \left( [ \ \underline{\theta} + (i-1) w, \ \underline{\theta} + i w \ ] \right),
\]
and we have 
\[
\begin{array}{lcl}
d(w \sum_{i=1}^{W} P^*(\Theta_W^{i})) = \mathcal{O} (1/W)
& \implies & \mathbf{A}^p_{\mathfrak{U}_W} = 1 - \mathcal{O} (1/W)
\end{array}
\]
Therefore the lower bound for the acceptance probability $\mathbf{A}^p_{\mathfrak{U}_W}$ 
of MRS approaches $1$ no slower than linearly with the refinement of 
$\mathbf{\Theta}$ by $\mathfrak{U}_W$.  
Note that this should hold for a general nonuniform partition with $w$ 
replaced by the mesh.$\square$

%% file: pqfRSMH.tex
\begin{picture}(0,0)%
\includegraphics{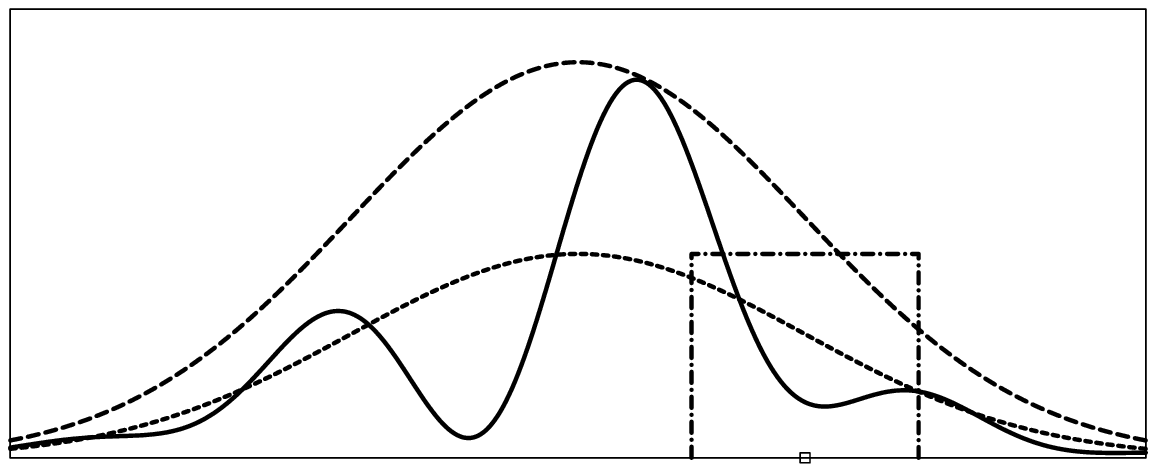}%
\end{picture}%
\begingroup
\setlength{\unitlength}{0.0200bp}%
\begin{picture}(18000,7560)(0,0)%
\put(8510,4866){\makebox(0,0)[l]{\strut{}$p^*$}}%
\put(6384,3340){\makebox(0,0)[l]{\strut{}$q^*_I$}}%
\put(13088,3956){\makebox(0,0)[l]{\strut{}$q^*_L$}}%
\put(7692,3780){\makebox(0,0)[l]{\strut{}$q^*$}}%
\put(11616,5101){\makebox(0,0)[l]{\strut{}$f_q$}}%
\end{picture}%
\endgroup
 

%% file: Gauss5_compare_PSzVar.tex
\begin{picture}(0,0)%
\includegraphics{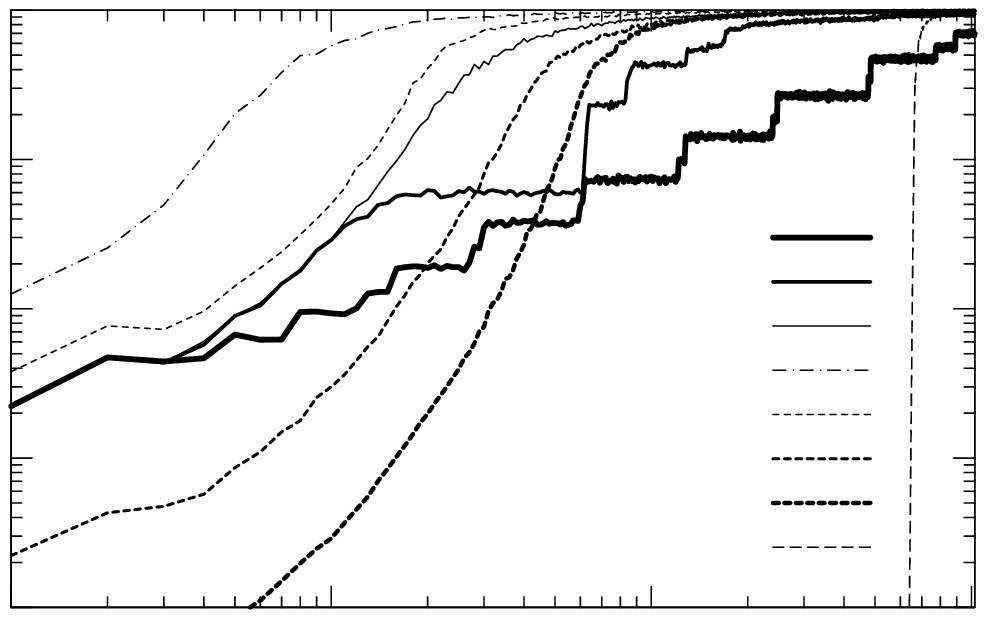}%
\end{picture}%
\begingroup
\setlength{\unitlength}{0.0200bp}%
\begin{picture}(18000,10800)(0,0)%
\put(3025,1650){\makebox(0,0)[r]{\strut{} 0.0001}}%
\put(3025,3800){\makebox(0,0)[r]{\strut{} 0.001}}%
\put(3025,5950){\makebox(0,0)[r]{\strut{} 0.01}}%
\put(3025,8100){\makebox(0,0)[r]{\strut{} 0.1}}%
\put(3025,10250){\makebox(0,0)[r]{\strut{} 1}}%
\put(3300,1100){\makebox(0,0){\strut{} 1}}%
\put(7909,1100){\makebox(0,0){\strut{} 10}}%
\put(12518,1100){\makebox(0,0){\strut{} 100}}%
\put(17128,1100){\makebox(0,0){\strut{} 1000}}%
\put(550,5950){\rotatebox{90}{\makebox(0,0){\strut{}Acceptance Probability ($\mathbf{A}^{g_5}_{\mathfrak{T}_{\alpha}}$)}}}%
\put(10237,275){\makebox(0,0){\strut{}Partition size ($|\mathfrak{T}_{\alpha}|$)}}%
\put(13996,6976){\makebox(0,0)[r]{\strut{}$g_5,\mathfrak{U}_{W}$}}%
\put(13996,6339){\makebox(0,0)[r]{\strut{}$g_5,\mathfrak{R}_{\alpha}$}}%
\put(13996,5702){\makebox(0,0)[r]{\strut{}$g_5,\mathfrak{V}_{\alpha}$}}%
\put(13996,5065){\makebox(0,0)[r]{\strut{}$g_1,\mathfrak{V}_{\alpha}$}}%
\put(13996,4428){\makebox(0,0)[r]{\strut{}$g_2,\mathfrak{V}_{\alpha}$}}%
\put(13996,3791){\makebox(0,0)[r]{\strut{}$g_5^{\prime},\mathfrak{V}_{\alpha}$}}%
\put(13996,3154){\makebox(0,0)[r]{\strut{}$g_5^{\prime \prime},\mathfrak{V}_{\alpha}$}}%
\put(13996,2517){\makebox(0,0)[r]{\strut{}$\widehat{g}_5,\mathfrak{V}_{\alpha}$}}%
\end{picture}%
\endgroup
 

%% file: Levy_compare_Temp.tex
\begin{picture}(0,0)%
\includegraphics{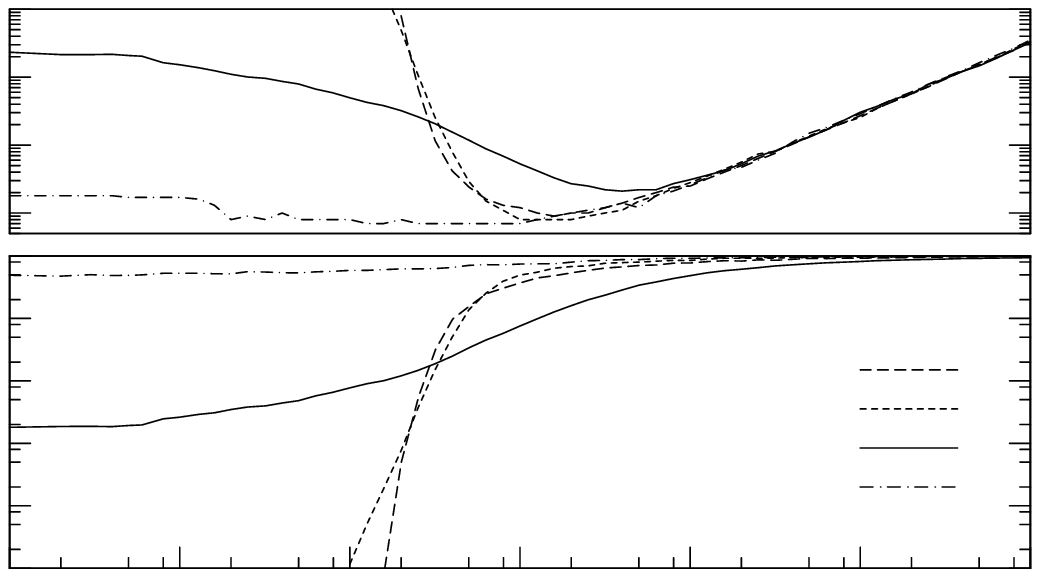}%
\end{picture}%
\begingroup
\setlength{\unitlength}{0.0200bp}%
\begin{picture}(18000,10800)(0,0)%
\put(1375,2200){\makebox(0,0)[r]{\strut{} 1e-05}}%
\put(1375,3099){\makebox(0,0)[r]{\strut{} 0.0001}}%
\put(1375,3998){\makebox(0,0)[r]{\strut{} 0.001}}%
\put(1375,4898){\makebox(0,0)[r]{\strut{} 0.01}}%
\put(1375,5797){\makebox(0,0)[r]{\strut{} 0.1}}%
\put(1375,6696){\makebox(0,0)[r]{\strut{} 1}}%
\put(1650,1650){\makebox(0,0){\strut{} 1}}%
\put(4100,1650){\makebox(0,0){\strut{} 10}}%
\put(6550,1650){\makebox(0,0){\strut{} 100}}%
\put(9000,1650){\makebox(0,0){\strut{} 1000}}%
\put(11450,1650){\makebox(0,0){\strut{} 10000}}%
\put(13900,1650){\makebox(0,0){\strut{} 100000}}%
\put(16350,1650){\makebox(0,0){\strut{} 1e+06}}%
\put(16899,4448){\rotatebox{90}{\makebox(0,0){\strut{}Accep.~Prob.~($\mathbf{A}^{l_T}_{\mathfrak{V}_{\alpha}}$)}}}%
\put(9000,825){\makebox(0,0){\strut{}Partition size ($|\mathfrak{V}_{\alpha}|$)}}%
\put(13625,5056){\makebox(0,0)[r]{\strut{}$l_1$}}%
\put(13625,4494){\makebox(0,0)[r]{\strut{}$l_4$}}%
\put(13625,3932){\makebox(0,0)[r]{\strut{}$l_{40}$}}%
\put(13625,3370){\makebox(0,0)[r]{\strut{}$l_{400}$}}%
\put(1375,7315){\makebox(0,0)[r]{\strut{} 0.1}}%
\put(1375,8293){\makebox(0,0)[r]{\strut{} 1}}%
\put(1375,9272){\makebox(0,0)[r]{\strut{} 10}}%
\put(1375,10250){\makebox(0,0)[r]{\strut{} 100}}%
\put(16899,8635){\rotatebox{90}{\makebox(0,0){\strut{}CPU Seconds}}}%
\end{picture}%
\endgroup
 

%% file: mix2norm_metro1.tex
\begin{picture}(0,0)%
\includegraphics{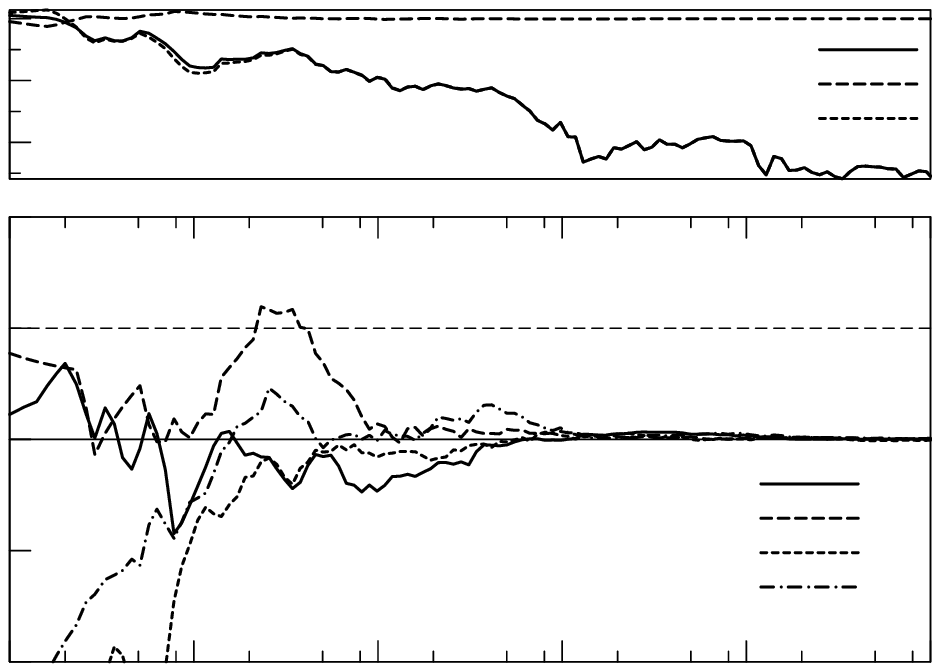}%
\end{picture}%
\begingroup
\setlength{\unitlength}{0.0200bp}%
\begin{picture}(18000,10800)(0,0)%
\put(1650,550){\makebox(0,0)[r]{\strut{}-1}}%
\put(1650,2151){\makebox(0,0)[r]{\strut{}-0.5}}%
\put(1650,3753){\makebox(0,0)[r]{\strut{} 0}}%
\put(1650,5354){\makebox(0,0)[r]{\strut{} 0.5}}%
\put(1650,6955){\makebox(0,0)[r]{\strut{} 1}}%
\put(1925,0){\makebox(0,0){\strut{} 10}}%
\put(4577,0){\makebox(0,0){\strut{} 100}}%
\put(7229,0){\makebox(0,0){\strut{} 1000}}%
\put(9881,0){\makebox(0,0){\strut{} 10000}}%
\put(12533,0){\makebox(0,0){\strut{} 100000}}%
\put(15185,0){\makebox(0,0){\strut{} 1e+06}}%
\put(15734,3752){\rotatebox{90}{\makebox(0,0){\strut{}Running Mean of $x_1$}}}%
\put(12468,3112){\makebox(0,0)[r]{\strut{}chain 1}}%
\put(12468,2617){\makebox(0,0)[r]{\strut{}chain 2}}%
\put(12468,2122){\makebox(0,0)[r]{\strut{}chain 3}}%
\put(12468,1627){\makebox(0,0)[r]{\strut{}chain 4}}%
\put(1650,8030){\makebox(0,0)[r]{\strut{} 0.0001}}%
\put(1650,8920){\makebox(0,0)[r]{\strut{} 0.01}}%
\put(1650,9810){\makebox(0,0)[r]{\strut{} 1}}%
\put(15734,8720){\rotatebox{90}{\makebox(0,0){\strut{}Variance of $x_1$}}}%
\put(13313,9365){\makebox(0,0)[r]{\strut{}B}}%
\put(13313,8870){\makebox(0,0)[r]{\strut{}W}}%
\put(13313,8375){\makebox(0,0)[r]{\strut{}B/W}}%
\end{picture}%
\endgroup
 

%% file: mrsmetrocompare.tex
\begin{picture}(0,0)%
\includegraphics{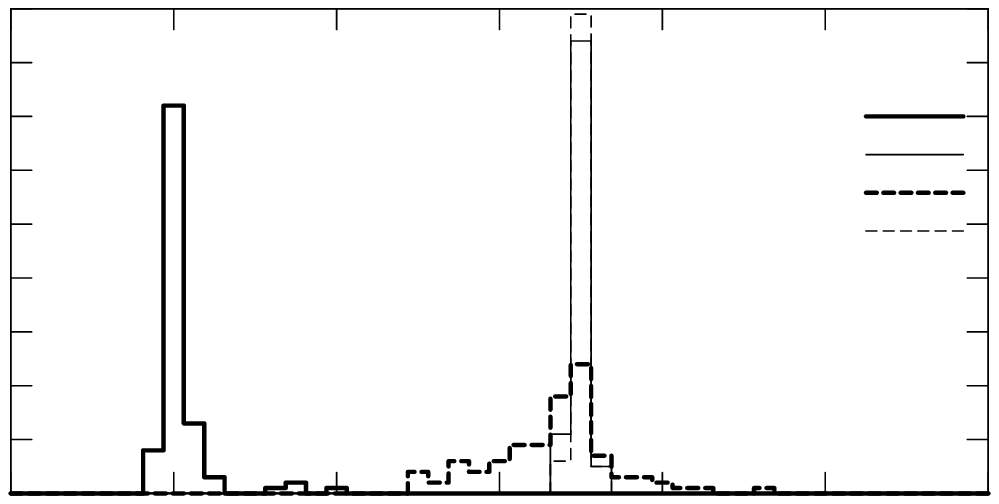}%
\end{picture}%
\begingroup
\setlength{\unitlength}{0.0200bp}%
\begin{picture}(17100,9180)(0,0)%
\put(1925,1650){\makebox(0,0)[r]{\strut{} 0}}%
\put(1925,2426){\makebox(0,0)[r]{\strut{} 10}}%
\put(1925,3201){\makebox(0,0)[r]{\strut{} 20}}%
\put(1925,3977){\makebox(0,0)[r]{\strut{} 30}}%
\put(1925,4752){\makebox(0,0)[r]{\strut{} 40}}%
\put(1925,5528){\makebox(0,0)[r]{\strut{} 50}}%
\put(1925,6303){\makebox(0,0)[r]{\strut{} 60}}%
\put(1925,7079){\makebox(0,0)[r]{\strut{} 70}}%
\put(1925,7854){\makebox(0,0)[r]{\strut{} 80}}%
\put(1925,8630){\makebox(0,0)[r]{\strut{} 90}}%
\put(2200,1100){\makebox(0,0){\strut{}-0.2}}%
\put(4546,1100){\makebox(0,0){\strut{} 0}}%
\put(6892,1100){\makebox(0,0){\strut{} 0.2}}%
\put(9238,1100){\makebox(0,0){\strut{} 0.4}}%
\put(11583,1100){\makebox(0,0){\strut{} 0.6}}%
\put(13929,1100){\makebox(0,0){\strut{} 0.8}}%
\put(16275,1100){\makebox(0,0){\strut{} 1}}%
\put(550,5140){\rotatebox{90}{\makebox(0,0){\strut{}Replicates}}}%
\put(9237,275){\makebox(0,0){\strut{}Mean $x_1$}}%
\put(14241,7079){\makebox(0,0)[r]{\strut{}LMHS(0.31)}}%
\put(14241,6529){\makebox(0,0)[r]{\strut{}MRS (0.15)}}%
\put(14241,5979){\makebox(0,0)[r]{\strut{}LMHS(3.36)}}%
\put(14241,5429){\makebox(0,0)[r]{\strut{}MRS (0.15)}}%
\end{picture}%
\endgroup
 

%% file: RSMISIS.tex
\begin{picture}(0,0)%
\includegraphics{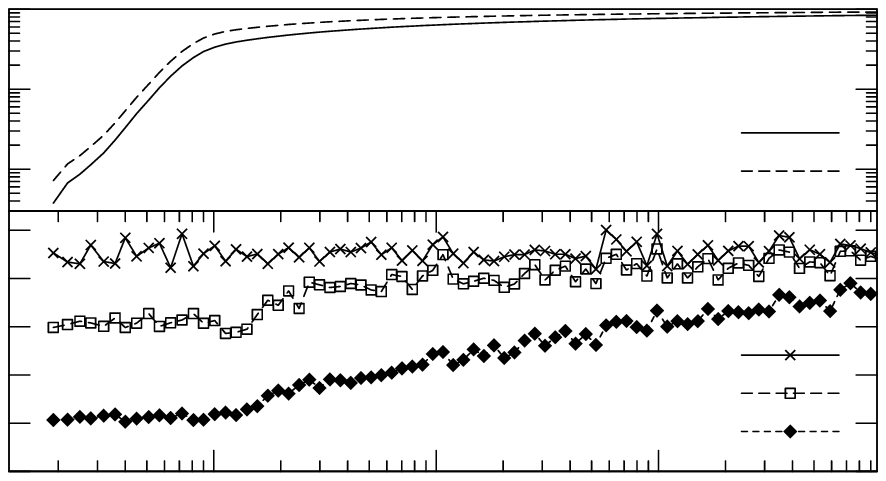}%
\end{picture}%
\begingroup
\setlength{\unitlength}{0.0200bp}%
\begin{picture}(18000,10800)(0,0)%
\put(2475,6002){\makebox(0,0)[r]{\strut{} 0.01}}%
\put(2475,7154){\makebox(0,0)[r]{\strut{} 0.1}}%
\put(2475,8306){\makebox(0,0)[r]{\strut{} 1}}%
\put(15799,6853){\rotatebox{90}{\makebox(0,0){\strut{}Accep.~Prob.}}}%
\put(13025,6525){\makebox(0,0)[r]{\strut{}MRS}}%
\put(13025,5975){\makebox(0,0)[r]{\strut{}IMHS}}%
\put(2475,1650){\makebox(0,0)[r]{\strut{} 0}}%
\put(2475,2344){\makebox(0,0)[r]{\strut{} 0.005}}%
\put(2475,3039){\makebox(0,0)[r]{\strut{} 0.01}}%
\put(2475,3733){\makebox(0,0)[r]{\strut{} 0.015}}%
\put(2475,4428){\makebox(0,0)[r]{\strut{} 0.02}}%
\put(2475,5122){\makebox(0,0)[r]{\strut{} 0.025}}%
\put(5699,1100){\makebox(0,0){\strut{} 100}}%
\put(8902,1100){\makebox(0,0){\strut{} 1000}}%
\put(12104,1100){\makebox(0,0){\strut{} 10000}}%
\put(15799,3525){\rotatebox{90}{\makebox(0,0){\strut{}MSE}}}%
\put(9000,275){\makebox(0,0){\strut{}Partition size ($|\mathfrak{V}_{\alpha}|$)}}%
\put(13025,3325){\makebox(0,0)[r]{\strut{}MRS}}%
\put(13025,2775){\makebox(0,0)[r]{\strut{}IMHS}}%
\put(13025,2225){\makebox(0,0)[r]{\strut{}IS}}%
\end{picture}%
\endgroup
 

%% file: Rosen_compare_Temp.tex
\begin{picture}(0,0)%
\includegraphics{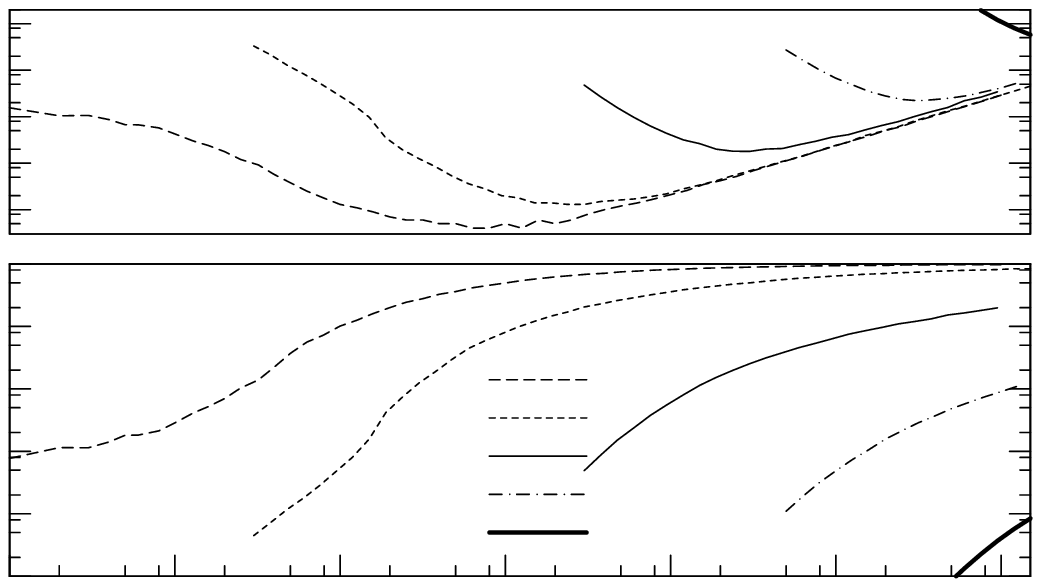}%
\end{picture}%
\begingroup
\setlength{\unitlength}{0.0200bp}%
\begin{picture}(18000,10800)(0,0)%
\put(1375,2200){\makebox(0,0)[r]{\strut{} 1e-05}}%
\put(1375,3099){\makebox(0,0)[r]{\strut{} 0.0001}}%
\put(1375,3998){\makebox(0,0)[r]{\strut{} 0.001}}%
\put(1375,4898){\makebox(0,0)[r]{\strut{} 0.01}}%
\put(1375,5797){\makebox(0,0)[r]{\strut{} 0.1}}%
\put(1375,6696){\makebox(0,0)[r]{\strut{} 1}}%
\put(1650,1650){\makebox(0,0){\strut{} 1}}%
\put(4030,1650){\makebox(0,0){\strut{} 10}}%
\put(6409,1650){\makebox(0,0){\strut{} 100}}%
\put(8789,1650){\makebox(0,0){\strut{} 1000}}%
\put(11169,1650){\makebox(0,0){\strut{} 10000}}%
\put(13549,1650){\makebox(0,0){\strut{} 100000}}%
\put(15928,1650){\makebox(0,0){\strut{} 1e+06}}%
\put(16899,4448){\rotatebox{90}{\makebox(0,0){\strut{}Accep.~Prob.~($\mathbf{A}^{r_D}_{\mathfrak{V}_{\alpha}}$)}}}%
\put(9000,825){\makebox(0,0){\strut{}Partition size ($|\mathfrak{V}_{\alpha}|$)}}%
\put(8284,5029){\makebox(0,0)[r]{\strut{}$r_2$}}%
\put(8284,4479){\makebox(0,0)[r]{\strut{}$r_3$}}%
\put(8284,3929){\makebox(0,0)[r]{\strut{}$r_5$}}%
\put(8284,3379){\makebox(0,0)[r]{\strut{}$r_7$}}%
\put(8284,2829){\makebox(0,0)[r]{\strut{}$r_9$}}%
\put(1375,7478){\makebox(0,0)[r]{\strut{} 0.1}}%
\put(1375,8148){\makebox(0,0)[r]{\strut{} 1}}%
\put(1375,8817){\makebox(0,0)[r]{\strut{} 10}}%
\put(1375,9487){\makebox(0,0)[r]{\strut{} 100}}%
\put(1375,10156){\makebox(0,0)[r]{\strut{} 1000}}%
\put(16899,8743){\rotatebox{90}{\makebox(0,0){\strut{}CPU Seconds}}}%
\end{picture}%
\endgroup
 

%% file: Wh1_compare_Dim_New.tex
\begin{picture}(0,0)%
\includegraphics{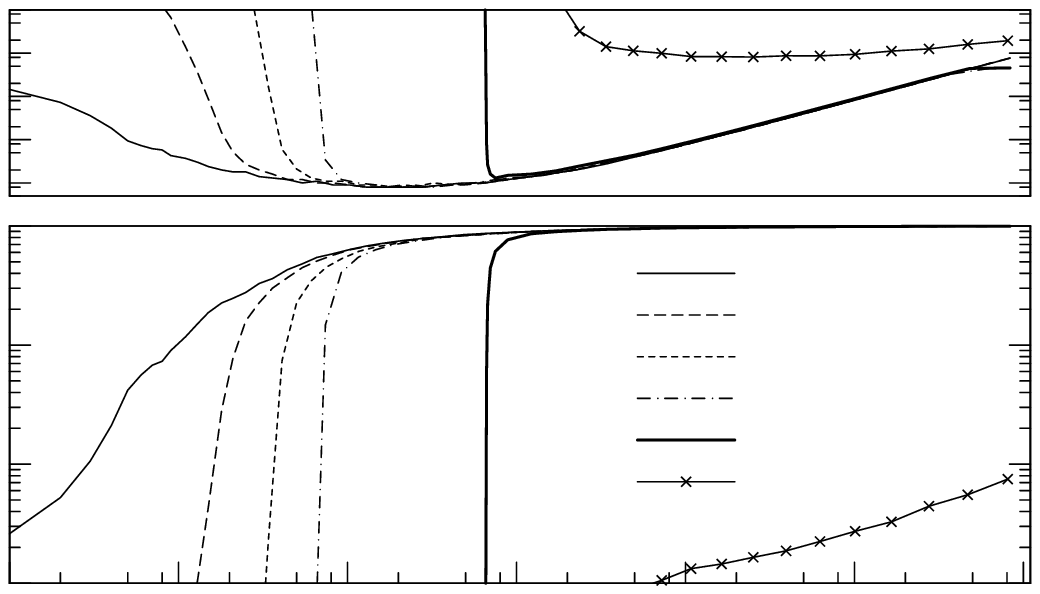}%
\end{picture}%
\begingroup
\setlength{\unitlength}{0.0200bp}%
\begin{picture}(18000,10800)(0,0)%
\put(1375,2200){\makebox(0,0)[r]{\strut{} 0.001}}%
\put(1375,3915){\makebox(0,0)[r]{\strut{} 0.01}}%
\put(1375,5629){\makebox(0,0)[r]{\strut{} 0.1}}%
\put(1375,7344){\makebox(0,0)[r]{\strut{} 1}}%
\put(1650,1650){\makebox(0,0){\strut{} 1}}%
\put(4083,1650){\makebox(0,0){\strut{} 10}}%
\put(6516,1650){\makebox(0,0){\strut{} 100}}%
\put(8950,1650){\makebox(0,0){\strut{} 1000}}%
\put(11383,1650){\makebox(0,0){\strut{} 10000}}%
\put(13816,1650){\makebox(0,0){\strut{} 100000}}%
\put(16249,1650){\makebox(0,0){\strut{} 1e+06}}%
\put(16899,4772){\rotatebox{90}{\makebox(0,0){\strut{}Accep.~Prob.~($\mathbf{A}^{h_r^D}_{\mathfrak{V}_{\alpha}}$)}}}%
\put(9000,825){\makebox(0,0){\strut{}Partition size ($|\mathfrak{V}_{\alpha}|$)}}%
\put(10417,6662){\makebox(0,0)[r]{\strut{}$h_0^2$}}%
\put(10417,6062){\makebox(0,0)[r]{\strut{}$h_2^2$}}%
\put(10417,5462){\makebox(0,0)[r]{\strut{}$h_5^2$}}%
\put(10417,4862){\makebox(0,0)[r]{\strut{}$h_{10}^2$}}%
\put(10417,4262){\makebox(0,0)[r]{\strut{}${\widehat{h}}_0^2$}}%
\put(10417,3662){\makebox(0,0)[r]{\strut{}$h_0^{10}$}}%
\put(1375,7964){\makebox(0,0)[r]{\strut{} 0.1}}%
\put(1375,8587){\makebox(0,0)[r]{\strut{} 1}}%
\put(1375,9210){\makebox(0,0)[r]{\strut{} 10}}%
\put(1375,9833){\makebox(0,0)[r]{\strut{} 100}}%
\put(1375,10456){\makebox(0,0)[r]{\strut{} 1000}}%
\put(16899,9116){\rotatebox{90}{\makebox(0,0){\strut{}CPU Seconds}}}%
\end{picture}%
\endgroup
 

%% file: DAG.pstex_t
\begin{picture}(0,0)%
\includegraphics{DAG.eps}%
\end{picture}%
\setlength{\unitlength}{3947sp}%
\begingroup\makeatletter\ifx\SetFigFont\undefined%
\gdef\SetFigFont#1#2#3#4#5{%
  \reset@font\fontsize{#1}{#2pt}%
  \fontfamily{#3}\fontseries{#4}\fontshape{#5}%
  \selectfont}%
\fi\endgroup%
\begin{picture}(6174,2625)(814,-2824)
\put(5551,-2761){\makebox(0,0)[lb]{\smash{{\SetFigFont{14}{16.8}{\rmdefault}{\mddefault}{\updefault}{\color[rgb]{0,0,0}$f_6=\cdot$}%
}}}}
\put(4801,-961){\makebox(0,0)[lb]{\smash{{\SetFigFont{14}{16.8}{\rmdefault}{\mddefault}{\updefault}{\color[rgb]{0,0,0}$s_5=\sin \left( \frac{x-3}{3} \right)$}%
}}}}
\put(5176,-2161){\makebox(0,0)[lb]{\smash{{\SetFigFont{14}{16.8}{\rmdefault}{\mddefault}{\updefault}{\color[rgb]{0,0,0}$s_6=x \sin \left( \frac{x-3}{3} \right)$}%
}}}}
\put(2101,-1411){\makebox(0,0)[lb]{\smash{{\SetFigFont{14}{16.8}{\rmdefault}{\mddefault}{\updefault}{\color[rgb]{0,0,0}$s_3=x-3$}%
}}}}
\put(976,-586){\makebox(0,0)[lb]{\smash{{\SetFigFont{14}{16.8}{\rmdefault}{\mddefault}{\updefault}{\color[rgb]{0,0,0}$s_2=3$}%
}}}}
\put(976,-2311){\makebox(0,0)[lb]{\smash{{\SetFigFont{14}{16.8}{\rmdefault}{\mddefault}{\updefault}{\color[rgb]{0,0,0}$s_1=x$}%
}}}}
\put(3301,-511){\makebox(0,0)[lb]{\smash{{\SetFigFont{14}{16.8}{\rmdefault}{\mddefault}{\updefault}{\color[rgb]{0,0,0}$s_4=\frac{x-3}{3}$}%
}}}}
\put(5026,-1486){\makebox(0,0)[lb]{\smash{{\SetFigFont{14}{16.8}{\rmdefault}{\mddefault}{\updefault}{\color[rgb]{0,0,0}$f_5=\sin$}%
}}}}
\put(3451,-1036){\makebox(0,0)[lb]{\smash{{\SetFigFont{14}{16.8}{\rmdefault}{\mddefault}{\updefault}{\color[rgb]{0,0,0}$f_4 = /$}%
}}}}
\put(2326,-1861){\makebox(0,0)[lb]{\smash{{\SetFigFont{14}{16.8}{\rmdefault}{\mddefault}{\updefault}{\color[rgb]{0,0,0}$f_3=-$}%
}}}}
\put(976,-2686){\makebox(0,0)[lb]{\smash{{\SetFigFont{14}{16.8}{\rmdefault}{\mddefault}{\updefault}{\color[rgb]{0,0,0}$f_1=s_1$}%
}}}}
\put(901,-961){\makebox(0,0)[lb]{\smash{{\SetFigFont{14}{16.8}{\rmdefault}{\mddefault}{\updefault}{\color[rgb]{0,0,0}$f_2=s_2$}%
}}}}
\end{picture}%